\documentclass[final,nopreprintline,3p]{elsarticle}
\usepackage{framed,multirow}

\usepackage{amssymb}
\usepackage{latexsym}

\usepackage{url}
\usepackage{xcolor}
\definecolor{newcolor}{rgb}{.8,.349,.1}

\usepackage{bm}
\usepackage{amsthm}
\usepackage{amsmath} 
\usepackage{amssymb}
\usepackage{subcaption}
\usepackage{tikz}
\usetikzlibrary{calc,patterns,matrix,shadings,fadings,shapes.arrows}
\usepackage{pgf-spectra}
\usepackage[Symbol]{upgreek}
\usepackage{enumitem}
\usepackage{tcolorbox}
\usepackage[ruled,vlined,linesnumbered,procnumbered]{algorithm2e}
\usepackage{multirow}
\usepackage{calc}
\usepackage{xcolor}
\usepackage{cancel}
\usepackage{booktabs}

\newcommand\vect[1]{\vec{#1}}
\usepackage{tikz-3dplot}
\tdplotsetmaincoords{60}{110}
\usepackage{mathtools}
\usepackage{adjustbox}
\usepackage{accents}

\usepackage{lmodern}
\usepackage[stretch=10]{microtype}

\usepackage{float}
\usepackage{etoolbox}

\usepackage{savesym}
\savesymbol{comment} 
\usepackage{comment}

\usepackage[outline]{contour}

\newsavebox{\foobox}

\newcommand\pgfmathsinandcos[3]{%
	\pgfmathsetmacro#1{sin(#3)}%
	\pgfmathsetmacro#2{cos(#3)}%
}
\newcommand\LongitudePlane[3][current plane]{%
	\pgfmathsinandcos\sinEl\cosEl{#2} 
	\pgfmathsinandcos\sint\cost{#3} 
	\tikzset{#1/.estyle={cm={\cost,\sint*\sinEl,0,\cosEl,(0,0)}}}
}
\newcommand\DrawLongitudeCircle[3][1]{
	\LongitudePlane{\angEl}{#2}
	\tikzset{current plane/.prefix style={scale=#1}}
	\pgfmathsetmacro\angVis{atan(sin(#2)*cos(\angEl)/sin(\angEl))} %
	\draw[current plane,#3] (0:1) arc (0:90:1);
}
\newcommand\DrawLongitudeCircleBoth[3][1]{
	\LongitudePlane{\angEl}{#2}
	\tikzset{current plane/.prefix style={scale=#1}}
	\pgfmathsetmacro\angVis{atan(sin(#2)*cos(\angEl)/sin(\angEl))} %
	\draw[current plane,#3] (0,-1) arc (-90:90:1);
}

\newcommand{\tikzcuboid}[4]{
	\begin{scope}[scale=#4]
		\foreach \x in {0,...,#1}
		{   
			\ifthenelse{\x=0 \OR \x=#1}{
				\colorlet{cx}{black}
			}{
				\colorlet{cx}{blue}
			}
			\draw[thick,cx] (\x ,0  ,#3 ) -- (\x ,#2 ,#3 );
			\draw[thick,cx] (\x ,#2 ,#3 ) -- (\x ,#2 ,0  );
		}
		\foreach \x in {0,...,#2}
		{   
			\ifthenelse{\x=0 \OR \x=#2}{
				\colorlet{cy}{black}
			}{
				\colorlet{cy}{green!50!black}
			}
			\draw[thick,cy] (#1 ,\x ,#3 ) -- (#1 ,\x ,0  );
			\draw[thick,cy] (0  ,\x ,#3 ) -- (#1 ,\x ,#3 );
		}
		\foreach \x in {0,...,#3}
		{   
			\ifthenelse{\x=0 \OR \x=#3}{
				\colorlet{cz}{black}
			}{
				\colorlet{cz}{red}
			}
			\draw[thick,cz] (#1 ,0  ,\x ) -- (#1 ,#2 ,\x );
			\draw[thick,cz] (0  ,#2 ,\x ) -- (#1 ,#2 ,\x );
		}
	\end{scope}
}

\newcommand{\tikzrectangle}[3]{
	\begin{scope}[#3]
		\fill[white,opacity=0.9] (0, 0) rectangle (#1, #2);
		\foreach \x in {0,...,#1} {
			\ifthenelse{\x=0 \OR \x=#1}{\colorlet{cx}{black}}{\colorlet{cx}{blue}}
			\draw[thick,cx] (\x, 0) -- (\x, #2);
		}
		\foreach \y in {0,...,#2} {
			\ifthenelse{\y=0 \OR \y=#2}{\colorlet{cy}{black}}{\colorlet{cy}{green!50!black}}
			\draw[thick,cy] (0, \y) -- (#1, \y);
		}
	\end{scope}
}

\begin{document}
	
\newcommand{\dir}{\alpha}
\newcommand{\ord}{\mu,\omega}

\begin{frontmatter}
	
	\title{Reduced-order modeling of neutron transport separated in axial and radial space by Proper Generalized Decomposition with applications to nuclear reactor physics}
	
	\fntext[fn1]{Present address: Naval Nuclear Laboratory, P.O. Box 1072, Schenectady, NY, 12301, USA.\\\\\raggedright\copyright 2024. This manuscript version is made available under the CC-BY-NC-ND 4.0 license \url{https://creativecommons.org/licenses/by-nc-nd/4.0/}}
	\cortext[cor2]{Corresponding author.}
	\author{Kurt A. {Dominesey}\fnref{fn1}}
	\ead{kurt.dominesey@unnpp.gov}
	\author{Wei {Ji}\corref{cor2}}
	\ead{jiw2@rpi.edu}
	
	\address{Department of Mechanical, Aerospace, \& Nuclear Engineering, Rensselaer Polytechnic Institute, Troy, NY 12180, USA}
	
	\begin{abstract}
		In this article, we demonstrate the novel use of Proper Generalized Decomposition (PGD) to separate the axial and, optionally, polar dimensions of neutron transport.
		Doing so, the resulting Reduced-Order Models (ROMs) 
		can exploit the fact that nuclear reactors tend to be tall, but geometrically simple, in the axial direction $z$, and so the 3D neutron flux distribution often admits a low-rank ``2D/1D'' approximation. 
		Through PGD, this approximation is computed by alternately solving 2D and 1D sub-models, like in existing 2D/1D models of reactor physics. 
		However, 
		the present methodology is more general in that the 
		decomposition is arbitrary-rank, rather than rank-one, and no simplifying approximations of the transverse leakage are made. 
		To begin, we derive two original models: that of axial PGD---which separates only $z$ and the sign of the polar angle $ \dir\in\{-1,+1\}$---and axial-polar PGD---which separates both $z$ and the full polar angle $\mu$ from the radial domain. Additionally, we grant that the energy dependence $E$ may be ascribed to either radial or axial modes, or both, bringing the total number of candidate 2D/1D ROMs to six.
		To assess performance, 
		these PGD ROMs are applied to two few-group benchmarks characteristic of Light Water Reactors. 
		Therein, we find both the axial and axial-polar ROMs are convergent and that the latter are often more economical than the former. 
		Ultimately, given the 
		popularity of 
		2D/1D methods in reactor physics, we expect a PGD ROM which achieves a similar effect, but perhaps with superior accuracy, a quicker runtime, and/or broader applicability, would be eminently useful, especially for full-core problems.
	\end{abstract}
	
\end{frontmatter}

\section{Motivation}
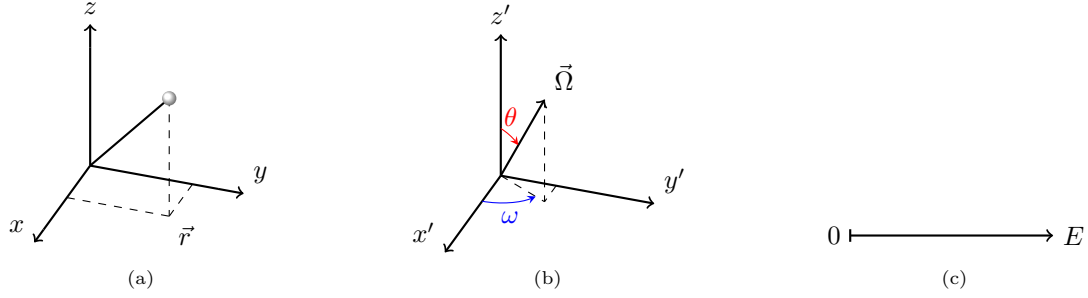
\begin{figure}[htbp]
	\begin{subfigure}[b]{0.32\textwidth}
		\centering
		\begin{tikzpicture}[scale=0.4*0.9,tdplot_main_coords]
		\def\Side{5}
		\def\X{2.5}
		\def\Y{4}
		\def\Z{5}
		
		\draw[thick] (0,0,0) -- (\Side,0,0);
		\draw[thick] (0,0,0) -- (0,\Side,0);
		\draw[thick] (0,0,0) -- (0,0,\Side);
		\draw[thick,->] (\Side,0,0) -- (\Side*1.2,0,0) node[anchor=south east]{$x$};
		\draw[thick,->] (0,\Side,0) -- (0,\Side*1.2,0) node[anchor=south west]{$y\hphantom{'}$};
		\draw[thick,->] (0,0,\Side) -- (0,0,\Side*1.2) node[anchor=south]{$z\vphantom{'}$};
		
		\coordinate (R) at (\X,\Y,\Z);
		\draw[dashed] (0,\Y,0) -- (\X,\Y,0);
		\draw[dashed] (\X,0,0) -- (\X,\Y,0);
		\draw[dashed] (R) -- (\X,\Y,0) node[anchor=north west] {$\vec{r}$};
		
		\draw[thick,-] (0,0,0) -- (R);
		\shade[ball color=white] (R) circle (7.5pt);
		\end{tikzpicture}
		\caption{}
	\end{subfigure}
	\begin{subfigure}[b]{0.32\textwidth}
		\centering
		\begin{tikzpicture}[scale=2*0.9,tdplot_main_coords]
		
		
		\def\rvec{1}
		\def\thetavec{30}
		\def\phivec{60}
		
		\def\nAzi{16}
		\def\nPol{4}
		\pgfmathsetmacro{\wedge}{360/\nAzi}
		\pgfmathsetmacro{\from}{\wedge/2}
		
		\coordinate (O) at (0,0,0);
		\draw[thick,->] (0,0,0) -- (1.2,0,0) node[anchor=south east]{$x'$};
		\draw[thick,->] (0,0,0) -- (0,1.2,0) node[anchor=south west]{$y'$};
		\draw[thick,->] (0,0,0) -- (0,0,1.2) node[anchor=south]{$z'$};
		
		\tdplotsetcoord{P}{\rvec}{\thetavec}{\phivec}
		\draw[thick,->] (O)  -- (P) node[anchor=south west] {$\vec{\Omega}$};
		\draw[dashed]   (O)  -- (Pxy);
		\draw[dashed]   (P)  -- (Pxy);
		\draw[dashed]   (Py) -- (Pxy);
		
		\tdplotdrawarc[-stealth,blue]{(O)}{0.4}{0}{\phivec}
		{anchor=north}{$\omega$}
		\tdplotsetthetaplanecoords{\phivec}
		\tdplotdrawarc[-stealth,red,tdplot_rotated_coords]{(0,0,0)}{0.4}{0}{\thetavec}
		{anchor=south west}{\hspace{-2mm}$\theta$}
		\end{tikzpicture}
		\caption{}
	\end{subfigure}
	\begin{subfigure}[b]{0.32\textwidth}
		\centering%
		\pgfdeclarehorizontalshading{rainbow}{100bp}{
			rgb(0bp)=(0,0,0);
			rgb(26bp)=(1,0,0);
			rgb(33bp)=(1,.5,0);
			rgb(40bp)=(1,1,0);
			rgb(47bp)=(0,1,0);
			rgb(54bp)=(0,1,1);
			rgb(61bp)=(0,0,1);
			rgb(68bp)=(1,0,1);
			rgb(75bp)=(.5,0,.5);
			rgb(100bp)=(0,0,0)}
		\begin{tikzpicture}[scale=0.21*0.9]
		\def\Ea{7}
		\def\Ec{-5}
		\def\tick{0.04*11}
		\fill[fill=white,shading=rainbow,draw=none,path fading=east] (\Ec,0) -- (\Ea+1.8,\Ea+1.8) -- (\Ea+1.8,0) -- cycle;
		\draw[thick,->] (\Ec,0) -- (\Ea+2.2,0) node[anchor=west] {$E$};
		\node[anchor=east] at (\Ec,0) {$0$};
		\draw[thick] (\Ec,-\tick) -- (\Ec,\tick);
		\end{tikzpicture}
		\caption{}
	\end{subfigure}
	\caption{The six-dimensional phase space of steady-state radiation transport consists of all possible spatial positions (a), angular trajectories (b), and energies (c).}
	\label{fig:phase-space}
\end{figure}

\newcommand{\twoD}{2\mathrm{D}}
\newcommand{\oneD}{1\mathrm{D}}
Computational simulation of radiation transport---as required to model, for instance, the neutronic inner workings of a nuclear reactor---is hampered by inherently high dimensionality, as depicted in Figure \ref{fig:phase-space}. That is, despite the governing equation being well-known and involving few assumptions,
calculating the neutron flux over all positions $\vect{r}$, angles $\vect{\Omega}$, and energies $E$ 
of interest in a nuclear reactor
can be intractable, even in two-dimensional (2D) space on a modern computing cluster. 
Given this, 
the computational burden may well become onerous or overwhelming in 3D geometry.

\label{sec:pgd_2d1d_motivation}
This is particularly true of
Light Water Reactors (LWRs),
which tend to be tall, though geometrically simple, in the axial, or $z$, dimension. To exploit this, as well as the short neutron mean-free-paths typical of LWRs, many ``2D/1D'' methods have been applied thereto---notably, those implemented within the CRX \cite{Cho2002}, DeCART \cite{Joo2004}, and MPACT \cite{Collins2016} simulators.
Each relies on various engineering assumptions---generally, some
homogenization of 2D unit-cells and a spatio-angular approximation of transverse leakage---to decouple the radial and axial dimensions, apart from an alternating iteration between 2D and 1D (sub-)problems.
While doing so renders 3D full-core transport tractable, the trade-off is some indeterminate error, numerical instability in certain circumstances,\footnote{Though Collins et al. \cite{Collins2016} show this can be rectified by a combination of Coarse Mesh Finite Difference (CMFD) acceleration and adaptive, energy dependent relaxation factors, as implemented in MPACT.} and limited applicability in the presence of severe flux gradients \cite{Jarrett2018}.

Inspired by, but methodologically distinct from, these 2D/1D approaches, we propose two analogous ROMs of neutron transport by Proper Generalized Decomposition (PGD) \cite{Ammar2006,Chinesta2010,Chinesta2014} as a means of decoupling radial and axial dimensions of neutron transport. First, we consider axial-polar PGD, where the neutron flux $\psi$ is separated as an expansion of $M$ radial and axial modes, $R_m$ and $Z_m^{(\mu)}$ respectively, 
\begin{equation}
\label{eq:axial-polar_pgd}
\psi(\vec{r},z,\mu,\omega,E)\approx\sum_{m=1}^M R_m(\vec{r},\omega,E) Z^{(\mu)}_m(z,\mu,E)\,,
\end{equation}
for radial position $\vec{r}\equiv(x,y)\in\mathcal{D}_2$, axial position $z\in\mathcal{D}_1\equiv[0,h]$, polar angle\footnote{More precisely, $\mu$ is the cosine of the polar angle $\theta$, as in $\mu\equiv\cos(\theta)$, but as $\theta$ is rarely used in this article, let us refer to $\mu$ as simply the polar angle.} $\mu\in[-1,1]$, azimuthal angle $\omega\in[0,2\pi)$, and energy $E\in[E_G,E_0]$. 
Second, we consider axial PGD, which is identical except that the polar angle is expressed as $\dir\mu$ for polar sign $\alpha=\pm 1$ and $\mu\in[0,1]$ so that $\alpha$ can be ascribed to the axial mode $Z_m$ and $\mu$ to the radial mode $R_m^{(\mu)}$,
\begin{equation}
\label{eq:axial_pgd}
\psi(\vec{r},z,\dir,\mu,\omega,E)=
\psi(\vec{r},z,\dir\mu,\omega,E)\approx\sum_{m=1}^M 
R^{(\mu)}_m
(\vec{r},\mu,\omega,E)
Z_m
(z,\dir,E)\,.
\end{equation}
The motivation for 
this change of coordinates is that in Equation \ref{eq:axial-polar_pgd} the sign of $\mu$ (positive or negative) dictates the axial streaming orientation (up or down) and therefore the boundary conditions on $Z_m$. 
Accordingly, $Z_m$ must be a function of the sign of the polar angle, if not its magnitude, which can be accomplished by respresenting these as separate variables $\alpha$ and $\mu$, as explained further in Section \ref{sec:sep_z_mu}.
Additionally, while Equations \ref{eq:axial-polar_pgd} and \ref{eq:axial_pgd} represent energy-dependent 2D/1D decompositions, in that both axial and radial modes are functions of energy, we also consider the alternative that $E$ could be solely ascribed to either of the two. As such, we present six candidate 2D/1D PGD ROMs: axial or axial-polar PGD, wherein energy is either shared between the axial and radial modes or else is separated out from one of the two.
These first two categories, axial- and axial-polar, are synonymously denoted as 2D($\mu$)/1D and 2D/1D($\mu$) where convenient.

Doing so, 
the radial and axial dimensions are again decoupled, notwithstanding an alternating iteration between 2D and 1D subproblems (outlined in Section \ref{sec:alg-pgd-2D1D}).
Crucially, however, this is accomplished with only a single assumption: that $\psi$ can be adequately approximated with a tractable number of modes $M$. This implies, also, that as $M$ grows large the model should converge to the true 3D solution.\footnote{This, of course, supposes that the approximation does converge, which is not theoretically guaranteed (see \cite{Cances2013}), but is often found in practice, including in the numerical results of Section \ref{sec:results-2D1D}.} 
This is not the case for 2D/1D methods generally, in that refining the 2D/1D mesh---as by adding additional 2D planes---does not typically guarantee convergence to the 3D flux (and may, in fact, render the scheme unstable without under-relaxation \cite{Kelley2015}), given other approximations remain.
Moreover, as these PGD ROMs do not rely on any properties of LWRs---save for the assumption of extruded geometry and the existence of a low-rank approximation---they are
immediately applicable to other reactors (graphite-moderated, fast, and so on)
and applications outside of reactor physics---perhaps including simulations of beamlines (for use in radiotherapy \cite{Han2011} or experimental physics) or solar atmospheric transport \cite{Davis2010}.

\section{Previous Work}
\label{sec:prev-work}
Mitigating the computation expense of 3D simulations is a longstanding issue in the reactor physics community. 
Arguably, the simplest approximation is to segment a 3D reactor into axial layers, in which the neutron source, geometry, and material properties are assumed to be axially constant; then, extrapolating the optical thickness of each layer to infinity, the original 3D problem simplifies to several 2D problems. 
Despite neglecting the effects of axial interfaces and boundaries completely, this geometric simplification is ubiquitous in deterministic cross section generation (that is, ``lattice physics'') and computational neutron transport in general \cite{Knott2010}.
Therefore, while not a 2D/1D method itself, the historical success of this approach reveals the essential feature of nuclear reactors that these methods exploit: namely, that the flux distribution of a 3D reactor can be adequately reconstructed from a relatively small set of 2D snapshots.

Some of the earliest efforts to exploit this property lie in ``flux synthesis'' methods \cite{Kaplan1962, Wachspress1962, Stacey1972, Kaplan1966, Alcouffe1970, Khromov1985} which decompose the solution into a sum of 2D and 1D modes, similar to those of Equations \ref{eq:axial_pgd} and \ref{eq:axial-polar_pgd}. Analysts would then provide one set of modes---either by physical intuition or based on previous reference solutions---and numerically compute the other: that is, compute the 1D modes given the 2D modes, or vice versa.\footnote{More specifically, analysts would provide 2D or 1D test and trial functions and apply a Petrov-Galerkin (or Galerkin if the test and trial functions are equivalent) projection to arrive at a 1D or 2D reduced-order model, respectively.} Doing so, the two can be ``synthesized'' into an approximate 3D solution.
This model was subsequently generalized to ``multichannel'' synthesis  wherein the decomposition is localized to particular radial zones (``channels'') as well as individual axial layers \cite{Wachspress1962, Wachspress1968,Stacey1968}.
Further,
some iterative versions of synthesis \cite{Wachspress1962,Slesarev1974}
alternated between solving for 1D and 2D modes until convergence, rendering the initial guess largely arbitrary. 
Despite these initial successes, however, synthesis methods 
are not widely used today
and important theoretical aspects---for instance, how to select test and trial functions---were not conclusively resolved.
That said, as flux synthesis is essentially a form of projection-based model order reduction---akin to Proper Orthogonal Decomposition (POD) \cite{Berkooz1993,Chatterjee2000} and PGD---much
relevant 
research has continued regardless, albeit largely 
conducted by other communities and applied to different ends. 

Modern 2D/1D models, meanwhile, 
proceed from 
the ``2D/1D fusion'' method of Cho et al. \cite{Cho2002} implemented in CRX. The defining feature is that the 3D geometry is divided into 1D layers and 2D zones (typically an assembly or pin-cell), wherein the flux is assumed to be exactly separable into a 2D and 1D component, both of which are computed numerically, iterating between the two until convergence. In this sense, they resemble an iterative, rank-one, multi-channel synthesis method as in \cite{Wachspress1962,Slesarev1974} or---perhaps even more similarly---Naito, Maekawa, and Shibuya's ``leakage iterative method'' \cite{Naito1975}.
However, several further approximations are typically made for convenience: for instance, that the transverse leakage---that is, the axial or radial leakage, alternately---is isotropic and spatially constant and that the ``homogenized'' cross sections are isotropic \cite{Cho2002,Joo2004,Collins2016, Stimpson2017, Jarrett2018, Kelley2015}.  
In principle, the PGD models developed in this article resemble these conventional 2D/1D methods, but generalized to an arbitrary number of modes (rather than a single one) and without such physically-motivated assumptions. 
That said, as further explanation requires some knowledge of PGD itself,
a more detailed description of specific similarities and differences 
is offered in Section \ref{sec:pgd-vs-2D1D}.

Finally, several projection-based reduced-order models---such as those of POD \cite{Berkooz1993,Chatterjee2000}, PGD \cite{Ammar2006, Chinesta2011, Chinesta2014}, Dynamic Mode Decomposition (DMD) \cite{Schmid2010,Tu2014}, and Dynamic Low-Rank Approximation (DLRA) \cite{Koch2007,Lubich2013,Ceruti2022}---have been applied to radiation transport and diffusion \cite{Buchan2013,Buchan2015,Lorenzi2018,Sun2020,Choi2021,Hardy2019,McClarren2019,McClarren2022,Roberts2019,Peng2020,Peng2021b,Peng2022,Einkemmer2021}. 
Of these, POD is methodologically similar to the flux synthesis methods already discussed, 
except that 
instead of  
relying on the practitioner to manually select
expansion modes (basis vectors), 
these are extracted from a set of reference solutions (snapshots) by a truncated matrix decomposition, namely a Singular Value Decomposition (SVD) or an equivalent eigendecomposition.
Meanwhile, PGD attempts to improve upon POD by eliminating the need for reference solutions, instead computing each mode progressively (that is, one at a time, beginning with zero modes) by an iterative procedure which alternates between low-dimensional sub-problems (in the present case, 2D and 1D subproblems) similar to that of iterative or ``closed'' synthesis \cite{Khromov1973,Slesarev1974,Khromov1985}. 

As applied to reactor physics, 
PGD has been previously used to separate $x$ and $y$ (1D/1D) in both diffusion \cite{Alberti2020,Senecal2018b,Senecal2019} and transport \cite{Prince2019} and $x$, $y$, and $z$ (1D/1D/1D) in diffusion \cite{Senecal2018,Prince2020,Prince2020b}. 
Further PGD separations, such as in time \cite{Alberti2019}, angle \cite{Dominesey2018,Dominesey2019a}, or energy \cite{Dominesey2019b,Prince2020,Dominesey2022,Dominesey2022b,Dominesey2023}, have also been demonstrated, sometimes in addition to the previous spatial separations.
However, a 2D/1D PGD method---capable of 
separating the 3D physics of a reactor
while faithfully representing the
unstructured 2D geometry---has as yet only been presented in Dominesey's dissertation \cite{DomineseyThesis}.
Additionally,
these prior works do not 
consider the simultaneous separation of both spatial and angular variables or the prospect of assigning $E$ to either low-order sub-model (rather than assigning it to both or separating out $E$ entirely).
Therefore, the objective of this article is to develop these 2D/1D PGD models and to investigate the practical and theoretical aspects of these features on prototypic LWR benchmarks.

\section{Methodology}
\label{sec:methodology-2D1D}
For brevity, we begin with the fixed-source neutron transport equation with $L$\textsuperscript{th}-order anisotropic scattering
\begin{equation}
\newlength{\hsep}
\label{eq:fom-2D1D}
\begin{split}
\sqrt{1-\mu^2}\left(\cos(\omega)\frac{\partial}{\partial x}+\sin(\omega)\frac{\partial}{\partial y}\right)\psi(\vec{r},z,\vec{\Omega},E)
+\mu\frac{\partial}{\partial z}\psi(\vec{r},z,\vec{\Omega},E)
+\Sigma_t(\vec{r},z,E)\psi(\vec{r},z,\vec{\Omega},E)
\\
-\sum_{\ell=0}^L\frac{2\ell+1}{4\pi}\sum_{k=-\ell}^{\ell}Y_{\ell,k}(\vec{\Omega})\int_{E_G}^{E_0}\Sigma_{s,\ell}(\vec{r},z,E'\rightarrow E)\int_{\mathcal{S}^2}Y_{\ell,k}(\vect{\Omega}')\psi(\vec{r},z,\vec{\Omega}',E')d\Omega'dE'
=q(\vec{r},z,\vec{\Omega},E)
\end{split}
\end{equation}
where $\Sigma_t$ and $\Sigma_s$ are macroscopic total and scattering cross sections, $Y_{\ell,k}$ are real spherical harmonics, $q$ is a neutron source, and $\vec{\Omega}$ is a direction on the unit sphere $\mathcal{S}^2$ defined by angular coordinates $(\mu,\,\omega)$. 
The separated models, to be solved for radial and axial modes, are then derived by first expanding each coefficient, as well as source $q$, analogously to $\psi$ in Equation \ref{eq:axial-polar_pgd} or \ref{eq:axial_pgd}, thereby
rendering Equation \ref{eq:fom-2D1D} separable in the sense required by either axial-polar or axial PGD.

\subsection{Separation of Coefficients}
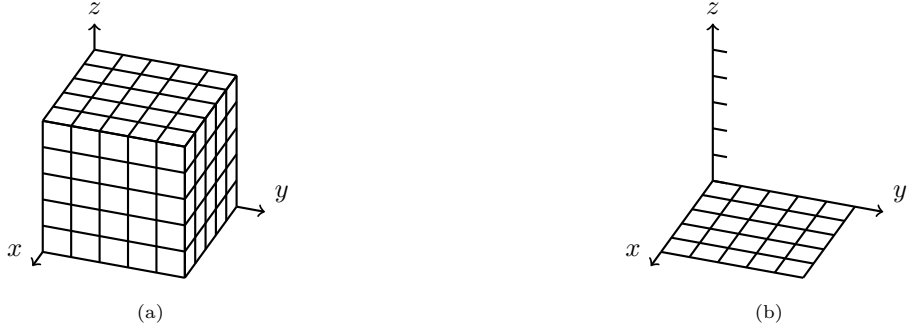
\begin{figure}[htpb]
	\centering
	\begin{subfigure}{0.49\textwidth}
		\centering
		\begin{tikzpicture}[scale=0.4,tdplot_main_coords]
		\def\Side{5}
		\def\Core{3}
		\def\Rod{4}
		\pgfmathsetmacro{\Refl}{\Side-\Core}
		\pgfmathsetmacro{\RodX}{\Side-\Rod}
		
		\node at (\Side,\Side,0) {};
		
		\draw[thick,->] (\Side,0,0) -- (\Side*1.2,0,0) node[anchor=south east]{$x$};
		\draw[thick,->] (0,\Side,0) -- (0,\Side*1.2,0) node[anchor=south west]{$y\hphantom{'}$};
		\draw[thick,->] (0,0,\Side) -- (0,0,\Side*1.2) node[anchor=south]{$z\vphantom{'}$};
		
		\foreach \i in {0,...,5} {
			\draw[thick] (\i,\Side,0) -- (\i,\Side,\Side);
		}
		\foreach \j in {0,...,5} {
			\draw[thick] (0,\Side,\j) -- (\Side,\Side,\j);
		}
		\foreach \i in {0,...,5} {
			\draw[thick] (\i,0,\Side) -- (\i,\Side,\Side);
		}
		\foreach \k in {0,...,5} {
			\draw[thick] (0,\k,\Side) -- (\Side,\k,\Side);
		}
		\foreach \j in {0,...,5} {
			\draw[thick] (\Side,0,\j) -- (\Side,\Side,\j);
		}
		\foreach \k in {0,...,5} {
			\draw[thick] (\Side,\k,0) -- (\Side,\k,\Side);
		}
		\end{tikzpicture}
		\caption{}
		\label{fig:cube-3D}
	\end{subfigure}
	\begin{subfigure}{0.49\textwidth}
		\centering
		\begin{tikzpicture}[scale=0.4,tdplot_main_coords]
		\def\Side{5}
		\def\Core{3}
		\def\Rod{4}
		\pgfmathsetmacro{\Refl}{\Side-\Core}
		\pgfmathsetmacro{\RodX}{\Side-\Rod}
		
		\node at (\Side,\Side,0) {};
		
		\draw[thick] (0,0,0) -- (0,0,\Side);
		\draw[thick,->] (\Side,0,0) -- (\Side*1.2,0,0) node[anchor=south east]{$x$};
		\draw[thick,->] (0,\Side,0) -- (0,\Side*1.2,0) node[anchor=south west]{$y\hphantom{'}$};
		\draw[thick,->] (0,0,\Side) -- (0,0,\Side*1.2) node[anchor=south]{$z\vphantom{'}$};
		
		\foreach \k in {0,...,5} {
			\draw [thick] (0,0,\k) -- (0,\Side*0.1,\k);
		}
		
		\foreach \i in {0,...,5} {
			\draw[thick] (\i,0,0) -- (\i,\Side,0);
		}
		\foreach \j in {0,...,5} {
			\draw[thick] (0,\j,0) -- (\Side,\j,0);
		}
		\end{tikzpicture}
		\caption{}
		\label{fig:cube-2D1D}
	\end{subfigure}
	\caption{Spatial decomposition of a 3D mesh (a) into 2D/1D meshes (b), as for Proper Generalized Decomposition.}
	\label{fig:cube}
\end{figure}

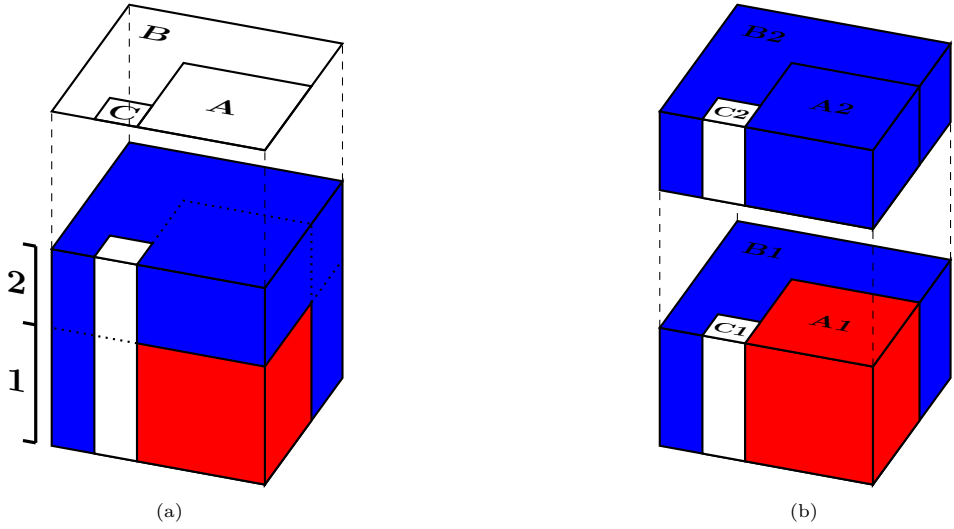
\begin{figure}
	\centering
	\tdplotsetmaincoords{60}{290}
	\begin{subfigure}{0.49\textwidth}
		\centering
		\begin{tikzpicture}[tdplot_main_coords,scale=0.6]
		\def\Side{5}
		\def\Core{3}
		\def\Rod{4}
		\pgfmathsetmacro{\Refl}{\Side-\Core}
		\pgfmathsetmacro{\RodX}{\Side-\Rod}
		
		\draw[thick,fill=blue] 
		(0,0,0) -- (\Side,0,0) -- (\Side,0,\Side) -- (0,0,\Side) -- cycle;
		\draw[thick,fill=blue] 
		(0,0,0) -- (0,\Side,0) -- (0,\Side,\Side) -- (0,0,\Side) -- cycle;
		\draw[thick,fill=blue] 
		(0,0,\Side) -- (0,\Side,\Side) -- (\Side,\Side,\Side) -- (\Side,0,\Side) -- cycle;
		
		\draw[thick,fill=red] 
		(0,0,0) -- (\Core,0,0) -- (\Core,0,\Core) -- (0,0,\Core) -- cycle;
		\draw[thick,fill=red] 
		(0,0,0) -- (0,\Core,0) -- (0,\Core,\Core) -- (0,0,\Core) -- cycle;
		;
		
		\draw[thick,fill=white] 
		(0,\Core,0) -- (0,\Rod,0) -- (0,\Rod,\Side) -- (0,\Core,\Side) -- cycle;
		\draw[thick,fill=white] 
		(0,\Core,\Side) -- (0,\Rod,\Side) -- (\RodX,\Rod,\Side) -- (\RodX,\Core,\Side) -- cycle;
		
		
		\draw[thick,dotted]
		(0,\Core,\Core) -- (0,\Side,\Core)
		(\Core,0,\Core) -- (\Side,0,\Core);
		
		\draw[thick,dotted] 
		(\RodX,\Core,\Side) -- (\Core,\Core,\Side) -- (\Core,0,\Side);
		\draw[thick,dotted] 
		(\Core,0,\Core) -- (\Core,0,\Side);
		
		\draw[dashed] 
		(0,0,\Side) -- +(0,0,3.5)
		(0,\Side,\Side) -- +(0,0,3.5)
		(\Side,0,\Side) -- +(0,0,3.5)
		(\Side,\Side,\Side) -- +(0,0,3.5);
		
		\begin{scope}[shift={(0,0,3.5)}]
		\draw[thick] 
		(0,0,\Side) -- (0,\Side,\Side) -- (\Side,\Side,\Side) -- (\Side,0,\Side) -- cycle;
		\node[canvas is xy plane at z=0,anchor=center,rotate=270] at (4,4,\Side) {\Large \textbf{B}};
		\draw[thick] 
		(0,\Core,\Side) -- (\Core,\Core,\Side) -- (\Core,0,\Side) -- (0,0,\Side) -- cycle;
		\node[canvas is xy plane at z=0,anchor=center,rotate=270] at (1.5,1.5,\Side) {\Large \textbf{A}};
		\draw[thick] 
		(0,\Core,\Side) -- (0,\Rod,\Side) -- (\RodX,\Rod,\Side) -- (\RodX,\Core,\Side) -- cycle;
		\node[canvas is xy plane at z=0,anchor=center,rotate=270] at (0.5,3.5,\Side) {\Large \textbf{C}};
		\end{scope}
		\begin{scope}[shift={(0,\Side*1.075,0)}]
		\draw[very thick] (0,0,0) -- (0,0,\Side);
		\draw[very thick]
		(0,0,0) -- +(0,0.3,0)
		(0,0,\Core) -- +(0,0.3,0)
		(0,0,\Side) -- +(0,0.3,0);
		\node[anchor=center,canvas is zy plane at x=0,rotate=270] at (0,0.45,1.5) {\Large\textbf{1}};
		\node[anchor=center,canvas is zy plane at x=0,rotate=270] at (0,0.45,4) {\Large\textbf{2}};
		\end{scope}
		\end{tikzpicture}
		\caption{}
		\label{fig:extruded-takeda}
	\end{subfigure}
	\hfill
	\begin{subfigure}{0.49\textwidth}
		\centering
		\begin{tikzpicture}[tdplot_main_coords,scale=0.6]
		\def\Side{5}
		\def\Core{3}
		\def\Rod{4}
		\pgfmathsetmacro{\Refl}{\Side-\Core}
		\pgfmathsetmacro{\RodX}{\Side-\Rod}
		
		\def\eZ{3}
		\draw[thick,fill=blue] 
		(0,0,0) -- (\Side,0,0) -- (\Side,0,\Core) -- (0,0,\Core) -- cycle;
		\draw[thick,fill=blue] 
		(0,0,0) -- (0,\Side,0) -- (0,\Side,\Core) -- (0,0,\Core) -- cycle;
		\draw[thick,fill=blue] 
		(0,0,\Core) -- (0,\Side,\Core) -- (\Side,\Side,\Core) -- (\Side,0,\Core) -- cycle;
		
		\draw[thick,fill=red] 
		(0,0,0) -- (\Core,0,0) -- (\Core,0,\Core) -- (0,0,\Core) -- cycle;
		\draw[thick,fill=red] 
		(0,0,0) -- (0,\Core,0) -- (0,\Core,\Core) -- (0,0,\Core) -- cycle;
		\draw[thick,fill=red] 
		(0,0,\Core) -- (0,\Core,\Core) -- (\Core,\Core,\Core) -- (\Core,0,\Core) -- cycle;
		
		\draw[thick,fill=white] 
		(0,\Core,0) -- (0,\Rod,0) -- (0,\Rod,\Core) -- (0,\Core,\Core) -- cycle;
		\draw[thick,fill=white] 
		(0,\Core,\Core) -- (0,\Rod,\Core) -- (\RodX,\Rod,\Core) -- (\RodX,\Core,\Core) -- cycle;
		
		\draw[dashed] 
		(0,0,\Core) -- +(0,0,3.5)
		(0,\Side,\Core) -- +(0,0,3.5)
		(\Side,0,\Core) -- +(0,0,3.5)
		(\Side,\Side,\Core) -- +(0,0,3.5);
		
		\node[canvas is xy plane at z=0,anchor=center,rotate=270] at (4,4,\Core) {\large\textbf{B1}};
		\node[canvas is xy plane at z=0,anchor=center,rotate=270] at (1.5,1.5,\Core) {\large\textbf{A1}};
		\node[canvas is xy plane at z=0,anchor=center,rotate=270] at (0.5,3.5,\Core) {\textbf{C1}};
		
		\begin{scope}[shift={(0,0,3.5)}]
		\path[thick] 
		(0,0,0) -- (0,\Side,0) -- (0,\Side,\Side) -- (0,0,\Side) -- cycle;
		
		\draw[thick,fill=blue] 
		(0,0,\Core) -- (\Side,0,\Core) -- (\Side,0,\Side) -- (0,0,\Side) -- cycle;
		\draw[thick,fill=blue] 
		(0,0,\Core) -- (0,\Side,\Core) -- (0,\Side,\Side) -- (0,0,\Side) -- cycle;
		
		\draw[thick,fill=blue] 
		(0,0,\Side) -- (0,\Side,\Side) -- (\Side,\Side,\Side) -- (\Side,0,\Side) -- cycle;
		\draw[thick,fill=blue] 
		(0,0,\Core) -- (\Core,0,\Core) -- (\Core,0,\Side) -- (0,0,\Side) -- cycle;
		\draw[thick,fill=blue] 
		(0,0,\Core) -- (0,\Core,\Core) -- (0,\Core,\Side) -- (0,0,\Side) -- cycle;
		\draw[thick,fill=blue] 
		(0,0,\Side) -- (0,\Core,\Side) -- (\Core,\Core,\Side) -- (\Core,0,\Side) -- cycle;
		
		\draw[thick,fill=white] 
		(0,\Core,\Core) -- (0,\Rod,\Core) -- (0,\Rod,\Side) -- (0,\Core,\Side) -- cycle;
		\draw[thick,fill=white] 
		(0,\Core,\Side) -- (0,\Rod,\Side) -- (\RodX,\Rod,\Side) -- (\RodX,\Core,\Side) -- cycle;
		
		\node[canvas is xy plane at z=0,anchor=center,rotate=270] at (4,4,\Side) {\large\textbf{B2}};
		\node[canvas is xy plane at z=0,anchor=center,rotate=270] at (1.5,1.5,\Side) {\large\textbf{A2}};
		\node[canvas is xy plane at z=0,anchor=center,rotate=270] at (0.5,3.5,\Side) {\textbf{C2}};
		\end{scope}
		\end{tikzpicture}
		\caption{}
	\end{subfigure}
	\caption{Extruded geometry of the Takeda Light Water Reactor (a);
		within each axial layer (1, 2) every radial area (A, B, C)
		corresponds to a homogeneous material region (b).}
	\label{fig:extruded-geom}
\end{figure}

To begin, let us assume the geoemtry of the 3D domain $\mathcal{D}^3\equiv\mathcal{D}^2\times\mathcal{D}^1$ admits a decomposition into {2D/1D}  sub-geometries, as in Figure \ref{fig:cube}---in other words, the 3D geometry is ``extrudable.''
Further, we grant the 2D and 1D domains $\mathcal{D}^2$ and $\mathcal{D}^1$ may be partitioned into $I$ radial areas $\mathcal{D}^2_i$ and $J$ axial layers $\mathcal{D}^1_j$ (respectively) such that each pair of indices ($i,j$) corresponds to a single, homogeneous material, exemplified in Figure \ref{fig:extruded-geom}.
Naturally, the cross sections may then be written
\begin{align}
\Sigma_{t}(\vec{r},z,E)&=\sum_{j=1}^J \sum_{i=1}^I\bm{1}_{\mathcal{D}_i^2}(\vec{r}) \bm{1}_{\mathcal{D}^1_j}(z) \Sigma_{t,i,j}(E) \\
\Sigma_{s,\ell}(\vec{r},z,E'\rightarrow E)&=\sum_{j=1}^J \sum_{i=1}^I  \bm{1}_{\mathcal{D}_i^2}(\vec{r})\bm{1}_{\mathcal{D}_j^1}(z) \Sigma_{s,\ell,i,j}(E'\rightarrow E)
\end{align}
where $\bm{1}$ is simply an indicator function
\begin{equation}
\bm{1}_{\mathcal{D}^2_i}(\vec{r})\equiv
\begin{cases}
1, & \vec{r} \in \mathcal{D}^2_i \\
0, & \mathrm{otherwise}
\end{cases}\,,
\qquad
\bm{1}_{\mathcal{D}^1_j}(z)\equiv
\begin{cases}
1, & z\in \mathcal{D}^1_j \\
0, & \mathrm{otherwise} \,.
\end{cases}
\end{equation}
Although these assumptions do limit the complexity of geometry which the ROM can represent,
many reactor physics problems can be expressed in this form, as 
is commonly exploited 
in other neutron transport software: for example, PROTEUS-MOC \cite{Jung2018} and MPACT \cite{Collins2016}.

Further, let us define the real spherical harmonics similarly to \cite{Hebert2009},
\begin{equation}
Y_{\ell,k}(\mu,\omega)\equiv
P^{|k|}_\ell(\mu) T_{k}(\omega)\,,
\end{equation}
where 
$P_\ell^{|k|}$ is the associated Legendre polynomial of degree $\ell$ and order $|k|$ multiplied by a normalization constant $\sqrt{\frac{(\ell-|k|)!}{(\ell+|k|)!}}$
and
\begin{equation}
T_{k}(\omega)\equiv
\begin{cases}
\sqrt{2}\sin(-k\omega)\,,
& k<0\,, \\
1
\,, & k=0\,, \\
\sqrt{2}\cos(k\omega)\,,
& k>0\,,
\end{cases}
\end{equation}
from which 
it is clear one can easily separate the polar and azimuthal dependencies, as required for axial-polar PGD. An additional step is needed for axial PGD, however, in that we have substituted $\dir\mu$ for $\mu$, but require $\dir$ and $\mu$ to be separated. This is overcome easily by noting
\begin{equation}
P_\ell^{|k|}(\dir\mu)=\dir^{\ell+k}P_\ell^{|k|}(\mu),\quad \forall\mu\in[0,1],\,\dir=\pm 1
\end{equation}
which effectively accounts for the parity of $P_{\ell}^{|k|}$ (either even or odd).

\subsection{Separation of Axial Space from Polar Angle}
\label{sec:sep_z_mu}
\begin{figure}[htbp]
	\centering
	\begin{subfigure}[t]{0.45\textwidth}{
			\begin{tikzpicture}
			\newcommand\xA{0}
			\newcommand\xB{6}
			\newcommand\muA{-1}
			\newcommand\muB{1}
			\draw[thin] (\xA,\muA) rectangle(\xB,\muB);
			\draw[thin,->] (\xA,\muA) -- (\xA,\muB+0.75) node[anchor=north east] {$\mu$};
			\foreach \mutick in {\muA,0,\muB}
			\draw (1pt,\mutick) -- (1pt,\mutick) node[anchor=east] {$\mutick$};
			\draw[thin,->] (\xA,\muA) -- (\xB+0.75,\muA) node[anchor=north east] {$z$};
			\draw (\xA,\muA) -- (\xA,\muA) node[anchor=north] {$0$};
			\draw (\xB,\muA) -- (\xB,\muA) node[anchor=north] {$h$};
			\draw[dashed] (\xA,0) -- (\xB,0);
			\draw[ultra thick,red] (\xA,0) -- (\xA,\muB);
			\draw[ultra thick,red] (\xB,0) -- (\xB,\muA);
			\draw[gray, thick, ->] (\xA+2.0,\muB/2) node[anchor=east,black] {$Z_m^{(\mu)}(z,\mu)$} --
			(\xA+2.66,\muB/2);
			\draw[gray, thick, ->] (\xA+2.0,\muB/2+0.3) --(\xA+3.0,\muB/2+0.3);
			\draw[gray, thick, ->] (\xA+2.0,\muB/2-0.3) --(\xA+2.33,\muB/2-0.3);
			\draw[gray, thick, ->] (\xB-2.0,\muA/2) node[anchor=west,black] {$Z_m^{(\mu)}(z,\mu)$} -- 
			(\xB-2.66,\muA/2);
			\draw[gray, thick, ->] (\xB-2.0,\muA/2+0.3) -- (\xB-2.33,\muA/2+0.3);
			\draw[gray, thick, ->] (\xB-2.0,\muA/2-0.3) -- (\xB-3.0,\muA/2-0.3);
			\draw[ultra thick,red] (\xA+0.4,\muB+0.75/2) -- (\xA+0.9,\muB+0.75/2)
			node[anchor=west] {Known boundary};
			\draw[thin] (\xA+4.3,\muB+0.75/2) -- (\xA+4.8,\muB+0.75/2)
			node[anchor=west] {Unknown boundary};
			\end{tikzpicture}
			\caption{}
			\label{fig:axial-polar-domain}
		}
	\end{subfigure}
	\hfill
	\begin{subfigure}[t]{0.45\textwidth}{
			\begin{tikzpicture}
			\newcommand\xA{0}
			\newcommand\xB{6}
			\newcommand\xM{(\xA/2+\xB/2}
			\newcommand\muA{0}
			\newcommand\muB{1}
			\draw[thin] (\xA,\muA) rectangle(\xB,\muB);
			\draw[thin,->] (\xA,\muA) -- (\xA,\muB+0.75) node[anchor=north east] {$\mu$};
			\foreach \mutick in {\muA,0,\muB}
			\draw (1pt,\mutick) -- (1pt,\mutick) node[anchor=east] {$\mutick$};
			\draw[thin,->] (\xA,\muA) -- (\xB+0.75,\muA) node[anchor=north east] {$z$};
			\draw (\xA,\muA) -- (\xA,\muA) node[anchor=north] {$0$};
			\draw (\xB,\muA) -- (\xB,\muA) node[anchor=north] {$h$};
			\draw[dashed] (\xA,0) -- (\xB,0);
			\draw[ultra thick,red,dashed] (\xA,0) -- (\xA,\muB);
			\draw[ultra thick,red,dashed] (\xB,0) -- (\xB,\muB);
			\draw[gray, thick, ->] (\xA+1.7,\muB/2) node[anchor=east,black] {$Z^{+}_m(z)$} --
			(\xA+2.4,\muB/2);
			\draw[gray, thick, ->] (\xB-1.7,\muB/2) node[anchor=west,black] {$Z^{-}_m(z)$} -- 
			(\xB-2.4,\muB/2);
			\draw[ultra thick,red,dashed] (\xA+0.4,\muB+0.75/2) -- (\xA+0.9,\muB+0.75/2)
			node[anchor=west] {Known to either $Z^{+}_m$ or $Z^{-}_m$};
			\end{tikzpicture}
			\caption{}
			\label{fig:axial-domain}
		}
	\end{subfigure}
	\caption{
		Axial-polar, or 1D($\mu$), (a) and axial, or 1D, (b) coordinate systems with gray arrows indicating particle direction.
	}
	\label{fig:pgd-domains}
\end{figure}
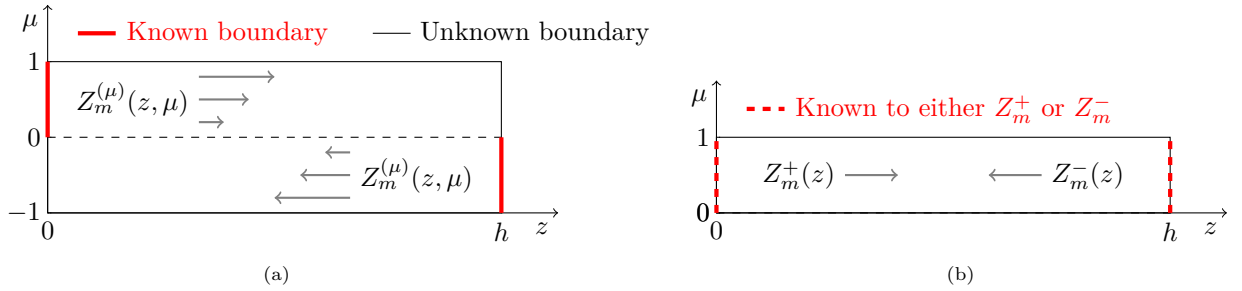

Although the axial position $z$ and polar angle $\mu$ appear to be separable in Equation \ref{eq:fom-2D1D},
attempting to separate these variables completely---as in, without introducing $\dir$---would render it impossible to impose the proper boundary conditions on the axial mode, $Z_m$.
The reason is that neutron transport is essentially an advective equation, and so the boundary conditions require some prescription of the incoming flux on the problem domain. However, at the axial (that is, top or bottom) boundaries, the polar angle dictates whether the flux is incoming (known) or outgoing (unknown); specifically, this is determined by the sign of $\mu$, either positive or negative. Plainly, without this information, there appears no means of computing $Z_m$. The straightforward resolution is to substitute 
$\mu$ with $\dir \mu$ 
where $\dir=\pm 1$ is the polar sign,
and to restrict $\mu$ from the original domain $[-1,+1]$ to $[0,1]$.
Likewise $\vec{\Omega}$, the direction specified by ($\mu$, $\omega$), is restricted to the upper hemisphere of the unit sphere $\mathcal{S}^2_+$.
Doing so, one arrives at a slightly modified equation of neutron transport
\begin{equation}
\setlength{\hsep}{2.25em}
\label{eq:fom-axial}
\begin{split}
\sqrt{1-\mu^2}\left(\cos(\omega)\frac{\partial}{\partial x}+\sin(\omega)\frac{\partial}{\partial y}\right)\psi(\vect{r},z,\dir,\vect{\Omega},E)
+\dir\mu\frac{\partial}{\partial z}\psi(\vect{r},z,\dir,\vect{\Omega},E)
+\Sigma_t(\vect{r},z,E)\psi(\vect{r},z,\dir,\vect{\Omega},E)
\\
-\sum_{\ell=0}^L\frac{2\ell+1}{4\pi} \sum_{k=-\ell}^{\ell} 
\dir^{\ell+k}Y_{\ell,k}(\vect{\Omega})
\int_{E_G}^{E_0} \Sigma_{s,\ell}(\vect{r},z,E'\rightarrow E)
\sum_{\dir'=\pm 1} \left(\dir'\right)^{\ell+k}\int_{\mathcal{S}^2_+} Y_{\ell,k}(\vect{\Omega}')\psi(\vect{r},z,\dir',\vect{\Omega}',E') d\Omega'dE' 
\\
=q(\vect{r},z,\dir,\vect{\Omega},E)
\end{split}
\end{equation}
which is amenable to separation by axial PGD. Meanwhile, this issue does not arise in axial-polar PGD since in that case the axial and polar dimensions are not separated.
These two different problem domains, including the requisite change of coordinates, are visualized in Figure  \ref{fig:pgd-domains}, where
\begin{equation}
Z_m(z,\dir)\equiv
\begin{cases}
Z^{+}_m(z), & \dir=+1 \,, \\
Z^{-}_m(z), & \dir=-1 \,.
\end{cases}
\end{equation}

Notably, this approach was first devised by Chinesta et al. \cite{Chinesta2015} and subsequently applied to neutron transport with collision and scattering in Dominesey and Ji \cite{Dominesey2019a}. Both studies also demonstrated agreement with analytical benchmarks.
As the present case of 2D($\mu$)/1D PGD is similar---except that azimuthal angle and radial space are no longer simplified away---these earlier works suggest the above strategy to be both straightforward and effective, as illustrated in the following sections.

\subsection{Separated Representation}
The derivation continues by substituting the separated representation of the flux---either that of axial-polar or axial PGD, Equation \ref{eq:axial-polar_pgd} or \ref{eq:axial_pgd}---and a similarly separated expansion of the source, either
\begin{equation}
\label{eq:q-axial-polar}
q(\vec{r},z,\mu,\omega,E)=\sum_{f=1}^{F}q_f^{\twoD}(\vec{r},\omega,E)q_f^{\oneD(\mu)}(z,\mu,E)
\end{equation}
or
\begin{equation}
\label{eq:q-axial}
q(\vec{r},z,\dir,\mu,\omega,E)= \sum_{f=1}^{F}q_f^{\twoD(\mu)}(\vec{r},\mu,\omega,E)q_f^{\oneD}(z,\dir,E)
\end{equation}
respectively, into Equation \ref{eq:fom-axial} or \ref{eq:fom-2D1D}. This decomposition of $q$ (where the rank $F$ is arbitrary) may be known \textit{a priori} or computed numerically, such as by the SVD (which can be truncated if retaining the full-rank decomposition is impractical). This, along with the separated cross sections and spherical harmonics, yields a statement of neutron transport which is compatible with 2D/1D PGD. Namely, in the axial-polar, 2D/1D($\mu$), case, substituting $\psi$ and $q$ as in Equations \ref{eq:axial-polar_pgd} and \ref{eq:q-axial-polar} into Equation \ref{eq:fom-2D1D} yields
\begin{equation}
\begin{split}
\sum_{m=1}^M\Bigg[\sqrt{1-\mu^2}Z^{(\mu)}_m(z,\mu,E)\times\left(\cos(\omega)\frac{\partial}{\partial x}+\sin(\omega)\frac{\partial}{\partial y}\right) R_m(\vec{r},\omega,E)
\\
+\mu\frac{\partial}{\partial z}Z^{(\mu)}_m(z,\mu,E)\times R_m(\vec{r},\omega,E)
\\
+\sum_{j=1}^J\sum_{i=1}^I
\Sigma_{t,i,j}(E)\times
\bm{1}_{\mathcal{D}^1_j}(z)Z_m^{(\mu)}(z,\mu,E)\times \bm{1}_{\mathcal{D}^2_i}(\vec{r})R_m(\vec{r},\omega,E) \\
-\sum_{\ell=0}^L\frac{2\ell+1}{4\pi}\sum_{k=-\ell}^{\ell}\sum_{j=1}^J\sum_{i=1}^I
\int_{E_G}^{E_0}dE'\,\Sigma_{s,\ell,i,j}(E'\rightarrow E)\\\times
\bm{1}_{\mathcal{D}^1_j}(z)P_{\ell}^{|k|}(\mu)\int_{-1}^{+1}P_{\ell}^{|k|}(\mu')Z_m^{(\mu)}(z,\mu',E')d\mu'\times \bm{1}_{\mathcal{D}^2_i}(\vec{r})T_{k}(\omega) \int_{0}^{2\pi}T_{k}(\omega')R_m(\vec{r},\omega',E')d\omega'\Bigg] \\
=\sum_{f=1}^F q^{\twoD}_f(\vec{r},\omega,E)q^{\oneD(\mu)}_f(z,\mu,E)
\end{split}
\end{equation}
where $\mu\in[-1,1]$
and likewise in the axial, 2D($\mu$)/1D, case, substituting Equations \ref{eq:axial_pgd} and \ref{eq:q-axial} into Equation \ref{eq:fom-axial} yields
\begin{equation}
\begin{split}
\sum_{m=1}^M\Bigg[Z_m(z,\dir,E)\times\sqrt{1-\mu^2}\left(\cos(\omega)\frac{\partial}{\partial x}+\sin(\omega)\frac{\partial}{\partial y}\right) R_m^{(\mu)}(\vec{r},\vec{\Omega},E)
\\
+\dir\frac{\partial}{\partial z}Z_m(z,\dir,E)\times\mu R_m^{(\mu)}(\vec{r},\vec{\Omega},E)
\\
+\sum_{j=1}^J\sum_{i=1}^I 
\Sigma_{t,i,j}(E)\times\bm{1}_{\mathcal{D}^1_j}(z)Z_m(z,\dir,E)\times\bm{1}_{\mathcal{D}^2_i}(\vec{r})R_m^{(\mu)}(\vec{r},\vec{\Omega},E)
\\
-\sum_{\ell=0}^L\frac{2\ell+1}{4\pi}\sum_{k=-\ell}^{\ell}\sum_{j=1}^J\sum_{i=1}^I \int_{E_G}^{E_0}dE'\,\Sigma_{s,\ell,i,j}(E'\rightarrow E)
\\\times
\bm{1}_{\mathcal{D}^1_j}(z)\dir^{\ell+k}\sum_{\dir'=\pm 1}\left(\dir'\right)^{\ell+k}Z_m(z,\dir',E')
\times
\bm{1}_{\mathcal{D}^2_i}(\vec{r})Y_{\ell,k}(\vec{\Omega})\int_{\mathcal{S}^2_+}Y_{\ell,k}(\vect{\Omega}')R_m^{(\mu)}(\vec{r},\vec{\Omega}',E')d\Omega'\Bigg]
\\
=\sum_{f=1}^F q^{\twoD(\mu)}_f(\vec{r},\vec{\Omega},E)q^{\oneD}_f(z,\dir,E)
\end{split}
\end{equation}
where $\mu\in[0,1]$.

Despite the change of notation,
these equations 
are 
equivalent to their original formulation in Equation \ref{eq:fom-2D1D}, save for the separated approximation of the solution. From this starting point, the PGD submodels---the 2D($\mu$), 2D, 1D, and 1D($\mu$) equations which are solved for individual radial or axial modes---are determined by a Galerkin projection, as detailed in Sections \ref{sec:theory-axial-polar-pgd} and \ref{sec:theory-axial-pgd}.
Then, having derived and discretized such submodels, PGD proceeds by 
way 
of an outer loop which repeatedly ``enriches'' the reduced basis (that is, computes the next term of the series); this entails, at each iteration, the resolution of
an inner loop which alternately solves each submodel in sequence until converging on a corresponding pair of modes. 
This procedure is outlined in Section \ref{sec:alg-pgd-2D1D}.

\subsection{Axial-Polar Decomposition, 2D/1D($\mu$)}
\label{sec:theory-axial-polar-pgd}
For simplicity, though without loss of generality, we here assume that energy $E$ has been discretized by the ubiquitous multigroup approximation.
Moreover, we detail only the most general case in full, where both submodels are energy-dependent;
simplifying either to the energy-independent case is straightforward, as demonstrated in Section \ref{sec:radial-axial-energy}. As will be seen, the resultant equations---in either the 2D/1D($\mu$) or 2D($\mu$)/1D case---can easily be recognized as typical statements of neutron transport in 1D and 2D geometry, notwithstanding certain peculiarities of the cross sections, source, and angular domains.
The details of the former---that is
the 2D and 1D($\mu$) submodels of the axial-polar PGD---are described in the following sections, while those of 2D($\mu$)/1D PGD are given in Section \ref{sec:theory-axial-pgd}.

\subsubsection{Radial Submodel, 2D}
To arrive at the radial, or 2D, submodel, 
consider the separable (full-order) model of neutron transport, Equation \ref{eq:fom-2D1D}, multiplied by test function $Z^{(\mu)}_{m^*}$ and integrated over the 1D spatial domain $\mathcal{D}^1$ and the polar axis $\mathcal{B}^1$, as in
\begin{equation}
\label{eq:proj-2D}
\int_{0}^{h}\int_{-1}^{+1} \bullet\ Z_{m^*,g}^{(\mu)}(z,\mu) d\mu dz\equiv
\left(Z^{(\mu)}_{m^*,g},\bullet\right)_{\mathcal{D}^1,\mathcal{B}^1}
\,.
\end{equation}
That is, we take the Galerkin projection of the full-order equation,
where   $Z^{(\mu)}_{m^*}$ is the axial-polar, or 1D($\mu$), mode indexed by $m^*$. In the Progressive PGD, only the equation corresponding to $m^*=M$ is solved, $M$ being progressively incremented to its final value, but $m^*$ is used here for generality. 
The cross sections in each radial area $\mathcal{D}^2_i$ can then be written as
\begin{align}
\widetilde{\Sigma}_{t,g,i}^{(m^*,m)}
&\equiv
\sum_{j=1}^J 
\left(Z_{m^*,g}^{(\mu)},Z_{m,g}^{(\mu)}\right)_{\mathcal{D}^1_j,\mathcal{B}^1}
\Sigma_{t,g,i,j}\,,
\\
\label{eq:sigma_sk}
\widetilde{\Sigma}_{s,k,g'\rightarrow g,i}^{(m^*,m)}
&\equiv 
\sum_{j=1}^J \sum_{\ell=|k|}^L\frac{2\ell+1}{2}\left(Z^{(\mu)}_{m^*,g},P_\ell^{|k|}\int_{-1}^{+1}P^{|k|}_\ell Z^{(\mu)}_{m,g'} d\mu'\right)_{\mathcal{D}^1_j,\mathcal{B}^1} \Sigma_{s,\ell,g'\rightarrow g,i,j}\,,
\end{align}
which amount to
flux-weighted sums over all layers $\mathcal{D}^1_j$.
Naturally, these constants appear in the separated expansion of the radial cross sections
\begin{align}
\widetilde{\Sigma}_{t,g}^{(m^*,m)}(\vec{r})&\equiv\sum_{i=1}^I \bm{1}_{\mathcal{D}^2_i}(\vec{r})\widetilde{\Sigma}_{t,g,i}^{(m^*,m)}
+\left(Z_{m^*,g}^{(\mu)},\mu\frac{\partial}{\partial z}Z^{(\mu)}_{m,g}\right)_{\mathcal{D}^1,\mathcal{B}^1}\,, \\
\widetilde{\Sigma}_{s,k,g'\rightarrow g}^{(m^*,m)}(\vec{r})&\equiv
\sum_{i=1}^I \bm{1}_{\mathcal{D}^2_i}(\vec{r})\widetilde{\Sigma}_{s,k,g'\rightarrow g,i}^{(m^*,m)}\,,
\end{align}
which,
in effect, simply serve to select the appropriate cross section $i$, with an additional term in the total cross section $\widetilde{\Sigma}_{t,g}^{(m^*,m)}$ to account for transverse (axial) leakage. For brevity, let us introduce the coplanar streaming coefficient
\begin{equation}
s_{\twoD,g}^{(m^*,m)}\equiv\left(Z^{(\mu)}_{m^*,g},\sqrt{1-\mu^2}Z^{(\mu)}_{m,g}\right)_{\mathcal{D}^1,\mathcal{B}^1}
\end{equation}
such that we can write our radial operator for group $g$---that is, one block row of the multigroup radial operator $\mathrm{B}^{(m^*,m)}_{\twoD}$---as
\begin{equation}
\begin{split}
\mathrm{B}_{\twoD,g}^{(m^*,m)}R_{m}&\equiv
s_{\twoD}^{(m^*,m)}
\left(
\cos(\omega)\frac{\partial}{\partial x}+
\sin(\omega)\frac{\partial}{\partial y}
\right)R_{m,g}(\vec{r},\omega)
+\widetilde{\Sigma}_{t,g}^{(m^*,m)}(\vec{r})R_{m,g}(\vec{r},\omega)\\
&\hphantom{\equiv}-\frac{1}{2\pi}\sum_{k=-L}^L 
T_{k}(\omega)
\sum_{g'=1}^G \widetilde{\Sigma}_{s,k,g'\rightarrow g}^{(m^*,m)}(\vec{r}) \int_0^{2\pi} T_{k}(\omega')R_{m,g'}(\vec{r},\omega') d\omega'\,.
\end{split}
\end{equation}
where we have rewritten the double summation $\sum_{\ell=0}^L\sum_{k=-\ell}^{\ell}$ as $\sum_{k=-L}^L\sum_{\ell=|k|}^L$ so that the summation over $\ell$ can be evaluated in Equation \ref{eq:sigma_sk}. 
Meanwhile, we introduce a term for the effective source
\begin{equation}
\widetilde{q}^{(m^*)}_{\twoD,g}(\vec{r},\omega)\equiv\sum_{f=1}^F\left(Z^{(\mu)}_{m^*,g},q_{f,g}^{\oneD{(\mu)}}\right)_{\mathcal{D}^1,\mathcal{B}^1}q_{f,g}^{\twoD}(\vec{r},\omega)\,.
\end{equation}
Doing so, we achieve a reduced, 2D version of the original transport equation which may be solved for radial mode $R_{m^*}$
\begin{equation}
\label{eq:prog-2D}
\begin{split}
\mathrm{B}_{\twoD}^{(m^*,m^*)}R_{m^*}&=\widetilde{q}_{\twoD}^{(m^*)}-\sum_{m\neq m^*}^M \mathrm{B}_{\twoD}^{(m^*,m)}R_{m} \\
&=\widetilde{r}_{\twoD}^{(m^*)}
\end{split}
\end{equation}
where $\widetilde{r}_{\twoD}^{(m^*)}$ is a residual term---that is, the source minus the contribution of other modes. For convenience, the coefficient of the streaming term can be canceled by amending our definitions
\begin{equation}
\begin{split}
\widetilde{\Sigma}_{t,g}^{(m^*,m^*)}(\vec{r})&\leftarrow\widetilde{\Sigma}_{t,g}^{(m^*,m^*)}
/s_{\twoD,g}^{(m^*,m^*)}\,, \\
\widetilde{\Sigma}_{s,k,g'\rightarrow g}^{(m^*,m^*)}(\vec{r})&\leftarrow\widetilde{\Sigma}_{s,k,g'\rightarrow g}^{(m^*,m^*)}(\vec{r})
/s_{\twoD,g}^{(m^*,m^*)}\,, \\
\widetilde{r}_{\twoD,g}^{(m^*)}(\vec{r},\omega)&\leftarrow\widetilde{r}_{\twoD,g}^{(m^*)}(\vec{r},\omega)
/s_{\twoD,g}^{(m^*,m^*)}\,,
\end{split}
\end{equation}
such that Equation \ref{eq:prog-2D} becomes identical to a conventional 2D fixed-source neutron transport problem save for three details: first, the cross sections and source are redefined; second, the polar angle is absent; and third, the atypical scattering kernel (for $L>0$). These peculiarities aside, this equation can be implemented and solved numerically via the same discretizations and algorithms used in the full-order model.

\subsubsection{Axial-Polar Submodel, 1D($\mu$)}
The derivation for the complementary submodel---the axial-polar, 1D($\mu$) equation---proceeds in much the same manner. One multiplies Equation \ref{eq:fom-2D1D} by $R_{m^*,g}$ and integrates over the 2D domain $\mathcal{D}^2$ and the unit circle $\mathcal{S}^1$, as in
\begin{equation}
\label{eq:proj-1D-mu}
\int_{\mathcal{D}^2} \int_{0}^{2\pi}
\bullet\ R_{m^*,g}(\vec{r},\omega)d\omega dr\equiv
\left(R_{m^*,g},\bullet\right)_{\mathcal{D}^2,\mathcal{S}^1}\,,
\end{equation}
rendering the cross sections for each layer $j=1\ldots J$ equal to
\begin{align}
\widetilde{\Sigma}_{t,g,j}^{(m^*,m)}&\equiv\sum_{i=1}^I \left(R_{m^*,g},R_{m,g}\right)_{\mathcal{D}^2_i,\mathcal{S}^1} \Sigma_{t,g,i,j}\,, \\
\label{eq:sigma_skl}
\widetilde{\Sigma}_{s,\ell,k,g'\rightarrow g,j}^{(m^*,m)}&\equiv\sum_{i=1}^I \frac{1}{2\pi}
\begin{cases}\displaystyle
\left(R_{m^*,g},\int_{0}^{2\pi}R_{m,g'}d\omega'\right)_{\mathcal{D}^2_i,\mathcal{S}^1} \Sigma_{s,\ell,g'\rightarrow g,i,j}\,,&k=0\,,\\\displaystyle
\sum_{p=\pm 1} \left(R_{m^*,g},T_{pk}\int_{0}^{2\pi}T_{pk}R_{m,g'}d\omega'\right)_{\mathcal{D}^2_i,\mathcal{S}^1} \Sigma_{s,\ell,g'\rightarrow g,i,j}\,,&k>0\,.\\
\end{cases}
\end{align}
These constants are incorporated into the axial, $z$-dependent cross sections
\begin{align}
\begin{split}
\widetilde{\Sigma}_{t,g}^{(m^*,m)}(z,\mu)&\equiv
\sum_{j=1}^J \bm{1}_{\mathcal{D}^1_j}(z)\widetilde{\Sigma}_{t,g,j}^{(m^*,m)}
\hphantom{\equiv}+\sqrt{1-\mu^2}\left(
R_{m^*,g},\cos(\omega)\frac{\partial}{\partial x}R_{m,g}+\sin(\omega)\frac{\partial}{\partial y}R_{m,g}\right)_{\mathcal{D}^2,\mathcal{S}^1}
\end{split}
\\
\widetilde{\Sigma}_{s,\ell,k,g'\rightarrow g}^{(m^*,m)}(z)&\equiv\sum_{j=1}^J \bm{1}_{\mathcal{D}^1_j}(z)\widetilde{\Sigma}_{s,\ell,k,g'\rightarrow g,j}^{(m^*,m)}\,,
\end{align}
where we find the effective total cross section $\widetilde{\Sigma}_{t,g}^{(m^*,m)}$ is necessarily anisotropic, owing to the fact that the transverse (radial) leakage term contained within depends on polar angle. This is unusual, but poses no great difficulty to implement. Moreover, while the memory burden of storing angularly-dependent cross sections is generally considered onerous, the same does not apply here, since the anisotropic coefficient $\sqrt{1-\mu^2}$ is analytical; as such, the additional storage requirement is solely that of the scalar inner product by which is it is multiplied. Then, given the coaxial streaming coefficient,
\begin{equation}
s_{\oneD(\mu),g}^{(m^*,m)}\equiv\left(R_{m^*,g},R_{m,g}\right)_{\mathcal{D}^2,\mathcal{S}^1}\,,
\end{equation}
the axial-polar operator for group $g$ can be written
\begin{equation}
\begin{split}
\mathrm{B}^{(m^*,m)}_{\oneD(\mu),g}Z^{(\mu)}_{m}&\equiv s_{\oneD(\mu),g}^{(m^*,m)}\mu\frac{\partial}{\partial z}Z_{m,g}^{(\mu)}(z,\mu)
+\widetilde{\Sigma}_{t,g}^{(m^*,m)}(z,\mu)Z_{m,g}^{(\mu)}(z,\mu) \\
&-\sum_{\ell=0}^L \frac{2\ell+1}{2}\sum_{k=0}^\ell P_{\ell}^{k}(\mu) \sum_{g'=1}^G \widetilde{\Sigma}_{s,\ell,k,g'\rightarrow g}^{(m^*,m)}(z)
\int_{-1}^{+1}P_{\ell}^{k}(\mu')Z_{m,g'}^{(\mu)}(z,\mu')d\mu'\,,
\end{split}
\end{equation}
where the summation $\sum_{k=-\ell}^\ell$ has been rewritten as $\sum_{k=0}^{\ell}$ since the terms for $k<0$ have been gathered into Equation \ref{eq:sigma_skl}.
Having done so, we now find a reduced, 1D neutron transport equation for the axial-polar mode $Z_{m^*}^{(\mu)}$
\begin{equation}
\label{eq:prog-1D-mu}
\begin{split}
\mathrm{B}_{\oneD(\mu)}^{(m^*,m^*)}Z^{(\mu)}_{m^*}&=\widetilde{q}_{\oneD(\mu)}^{(m^*)}-\sum_{m\neq m^*}^M \mathrm{B}_{\oneD(\mu)}^{(m^*,m)}Z^{(\mu)}_{m} \\
&=\widetilde{r}_{\oneD(\mu)}^{(m^*)}
\end{split}
\end{equation}
where the effective source is given by
\begin{equation}
\widetilde{q}^{(m^*)}_{\oneD(\mu),g}(z,\mu)\equiv\sum_{f=1}^F \left(R_{m^*,g},q^{\twoD}_{f,g}\right)_{\mathcal{D}^2,\mathcal{S}^1} q^{\oneD(\mu)}_{f,g}(z,\mu)\,.
\end{equation}
As before, the streaming coefficient can be canceled by substituting
\begin{equation}
\begin{split}
\widetilde{\Sigma}_{t,g}^{(m^*,m)}(z,\mu)&\leftarrow\widetilde{\Sigma}_{t,g}^{(m^*,m)}(z,\mu)/s_{\oneD(\mu),g}^{(m^*,m^*)}\,, \\
\widetilde{\Sigma}_{s,\ell,k,g'\rightarrow g}^{(m^*,m)}(z)&\leftarrow\widetilde{\Sigma}_{s,\ell,k,g'\rightarrow g}^{(m^*,m)}(z)/ s_{\oneD(\mu),g}^{(m^*,m^*)}\,, \\
\widetilde{r}^{(m^*)}_{\oneD(\mu),g}(z,\mu)&\leftarrow\widetilde{r}^{(m^*)}_{\oneD(\mu),g}(z,\mu)/s_{\oneD(\mu),g}^{(m^*,m^*)}\,,
\end{split}
\end{equation}
to render Equation \ref{eq:prog-1D-mu} a 1D fixed-source neutron transport problem which, besides the specification of cross sections and source, is unusual only in the presence of an anisotropic transverse leakage term and the inclusion of azimuthal components of the spherical harmonic moments other than $k=0$ in the scattering term. Specifically, the latter detail is unusual because the 1D solution is typically azimuthally constant by definition, such that these moments equal zero when integrated in $\omega$ over the unit circle $\mathcal{S}^1$.

\subsection{Axial Decomposition, 2D($\mu$)/1D}
\label{sec:theory-axial-pgd}
Turning our attention to the axial (rather than axial-polar), 2D($\mu$)/1D decomposition, 
we now start from Equation \ref{eq:fom-axial}, where, compared to Equation \ref{eq:fom-2D1D}, we have restricted the polar angle $\mu$ from $[-1,+1]$ to $[0,1]$ and introduced an additional coordinate for the polar sign $\dir\in\{-1,+1\}$ which represents the axial streaming orientation, up or down. As explained in Section \ref{sec:sep_z_mu}, this is necessary to separate $z$ and $\mu$ by PGD, and so is required if $\mu$ is to be assigned to the 2D submodel---as is the defining feature of the axial decomposition.
Apart from this distinction, the derivation proceeds similarly as in the preceding case of 2D/1D($\mu$) decomposition.

\subsubsection{Radial-Polar Submodel, 2D($\mu$)}
Specifically, to achieve the radial-polar, or 2D($\mu$), sub-model,
we multiply Equation \ref{eq:fom-axial} by $Z_{m^*,g}$, integrate over the 1D domain $\mathcal{D}^1$, and sum over both polar signs $\alpha\in\mathcal{S}^0\equiv\{-1,+1\}$ as in
\begin{equation}
\label{eq:proj-2D-mu}
\sum_{\dir=\pm 1}\int_{0}^h \bullet\,Z_{m^*,g}
(z,\dir)dz\equiv\left(Z_{m^*,g},\bullet\right)_{\mathcal{D}^1,\mathcal{S}^0}
\end{equation}
and subsequently define the axially-integrated, flux-weighted cross sections for each radial area $\mathcal{D}^2_i$
\begin{align}
\widetilde{\Sigma}_{t,g,i}^{(m^*,m)}&\equiv\sum_{j=1}^J \left(Z_{m^*,g},Z_{m,g}\right)_{\mathcal{D}^1_j,\mathcal{S}^0} \Sigma_{t,i,j,g}\,, \\
\widetilde{\Sigma}_{s,\ell,p,g'\rightarrow g,i}^{(m^*,m)}&\equiv\sum_{j=1}^J  
\frac{1}{2}\left(Z_{m^*,g},\dir^{p}\sum_{\dir'=\pm 1}\left(\dir'\right)^pZ_{m,g'}\right)_{\mathcal{D}^1_j,\mathcal{S}^0}
\Sigma_{s,\ell,g'\rightarrow g,i,j}\,,
\end{align}
where $p\equiv(\ell+k)\bmod2$ is the polar parity.
Further, by evaluating the summations over $\dir$ and $\dir'$, the latter inner product can be simplified as
\begin{equation}
\label{eq:parity-ip}
\left(Z_{m^*,g},\dir^{p}\sum_{\dir'=\pm 1}\left(\dir'\right)^{p}Z_{m,g'}\right)_{\mathcal{D}^1_j,\mathcal{S}^0}\gets
\left(Z_{m^*,g}^p,Z^p_{m,g'}\right)_{\mathcal{D}^1_j}
\end{equation}
where $Z_{m,g}^p$ is the even- or odd-parity axial flux moment,
\begin{equation}
Z^p_{m,g}(z)\equiv\begin{cases}
Z^{0}_{m,g}(z)\equiv Z^+_{m,g}(z)+Z^-_{m,g}(z)\,, & \ell+k\text{ is even}\,, \\
Z^{1}_{m,g}(z)\equiv Z^+_{m,g}(z)-Z^-_{m,g}(z)\,, & \ell+k\text{ is odd}\,.
\end{cases}
\end{equation}
These coefficients are then incorporated into the radial-polar cross sections
\begin{align}
\widetilde{\Sigma}_{t,g}^{(m^*,m)}(\vec{r},\mu)&\equiv\sum_{i=1}^I \bm{1}_{\mathcal{D}^2_i}(\vec{r})\widetilde{\Sigma}_{t,g,i}^{(m^*,m)}
+\mu\left(Z_{m^*,g},\dir\frac{\partial}{\partial z}Z_{m,g}\right)_{\mathcal{D}^1,\mathcal{S}^0}\,, \\
\widetilde{\Sigma}_{s,\ell,p,g'\rightarrow g}^{(m^*,m)}(\vec{r})&\equiv\sum_{i=1}^I \bm{1}_{\mathcal{D}^2_i}(\vec{r}) \widetilde{\Sigma}_{s,\ell,p,g'\rightarrow g,i}^{(m^*,m)}\,,
\end{align}
where, as in the axial-polar case, the total cross section $\widetilde{\Sigma}_{t,g}^{(m^*,m)}$ is a function of $\mu$ because the transverse (axial) leakage is dependent on polar angle.
Using these terms, plus a coplanar streaming coefficient, 
\begin{equation}
s^{(m^*,m)}_{\twoD(\mu),g}\equiv\left(Z_{m^*,g},Z_{m,g}\right)_{\mathcal{D}^1,\mathcal{S}^0}
\end{equation}
we define $\mathrm{B}^{(m^*,m)}_{\twoD(\mu),g}$, the radial-polar operator for group $g$
\begin{equation}
\begin{split}
\mathrm{B}^{(m^*,m)}_{\twoD(\mu),g} R^{(\mu)}_{m}&\equiv
s^{(m^*,m)}_{\twoD(\mu),g}
\sqrt{1-\mu^2}\left(\cos(\omega)\frac{\partial}{\partial x}+\sin(\omega)\frac{\partial}{\partial y}\right)
R_{m,g}^{(\mu)}(\vect{r},\vect{\Omega})
+\widetilde{\Sigma}_{t,g}^{(m^*,m)}(\vect{r},\mu)R_{m,g}^{(\mu)}(\vect{r},\vect{\Omega}) \\
&\hphantom{\equiv}-\sum_{\ell=0}^L \frac{2\ell+1}{2\pi}\sum_{k=-\ell}^{\ell}  Y_{\ell,k}(\vect{\Omega})
\sum_{g'=1}^G \widetilde{\Sigma}_{s,\ell,p,g'\rightarrow g}^{(m^*,m)}(\vec{r})
\int_{\mathcal{S}^2_+} Y_{\ell,k}(\vect{\Omega}')
R_{m,g'}^{(\mu)}
(\vect{r},\vect{\Omega}')d\Omega'
\end{split}
\end{equation}
which allows us to succinctly express our equation for $R^{(\mu)}_{m^*}$
\begin{equation}
\label{eq:prog-2D-mu}
\begin{split}
\mathrm{B}^{(m^*,m^*)}_{\twoD(\mu)} R^{(\mu)}_{m^*}&=\widetilde{q}^{(m^*)}_{\twoD(\mu),g}
-\sum_{m\neq m^*}\mathrm{B}^{(m^*,m)}_{\twoD(\mu)} R^{(\mu)}_{m} \\
&=\widetilde{r}^{(m^*)}_{\twoD(\mu)}
\end{split}
\end{equation}
where
\begin{equation}
\widetilde{q}^{(m^*)}_{\twoD(\mu),g}(\vect{r},\vect{\Omega})\equiv\sum_{f=1}^F
\left(Z_{m^*,g},q_{f,g}^{\oneD}\right)_{\mathcal{D}^1,\mathcal{S}^0} q_{f,g}^{\twoD(\mu)}(\vect{r},\vect{\Omega})\,.
\end{equation}
Again, to cancel the streaming coefficient, we can further redefine
\begin{equation}
\begin{split}
\widetilde{\Sigma}_{t,g}^{(m^*,m)}(\vec{r},\mu)&\leftarrow\widetilde{\Sigma}_{t,g}^{(m^*,m)}(\vec{r},\mu)
/s_{\twoD(\mu),g}^{(m^*,m^*)}\,, \\
\widetilde{\Sigma}_{s,\ell,p,g'\rightarrow g}^{(m^*,m)}(\vec{r})&\leftarrow\widetilde{\Sigma}_{s,\ell,p,g'\rightarrow g}^{(m^*,m)}(\vec{r})
/s_{\twoD(\mu),g}^{(m^*,m^*)}\,, \\
\widetilde{r}_{\twoD(\mu),g}^{(m^*)}(\vec{r},\vect{\Omega})&\leftarrow\widetilde{r}_{\twoD(\mu),g}^{(m^*)}(\vec{r},\vect{\Omega})
/s_{\twoD(\mu),g}^{(m^*,m^*)}\,,
\end{split}
\end{equation}
yielding an equation nearly identical with the usual statement of 2D fixed-source neutron transport. Specifically, it differs---beyond the redefinition of cross sections and source---in that the total cross section is anisotropic, and in that those spherical harmonics moments in the scattering term which are odd functions of $\mu$ (that is, where $\ell+k$ is odd and so $p=1$) do not integrate to zero. Of course, these moments vanish regardless if the corresponding cross sections becomes zero, such as when there is no net axial current, $Z^{1}_{m^*,g}=0$ or $Z^{1}_{m,g'}=0$, so that the inner product in Equation \ref{eq:parity-ip} vanishes. Further, while one may expect that solving for only half the range of polar angles---those in the interval [0,1]---is unusual, 
the case is in fact the opposite:
since the solution of any 2D problem will be symmetric with respect to $\mu$ (so long as the source is as well), there is rarely a need to solve for both (north and south) hemispheres. As such, this particular aspect of Equation \ref{eq:prog-2D-mu} should not pose any implementation difficulties.

\subsubsection{Axial Submodel, 1D}
Finding the corresponding axial equation, meanwhile, begins by multiplying Equation \ref{eq:fom-axial} by $R_{m^*,g}$ and integrating over the 2D domain $\mathcal{D}^2$ and the upper unit hemisphere $\mathcal{S}^2_+$ as in
\begin{equation}
\label{eq:proj-1D}
\int_{\mathcal{D}^2} \int_{\mathcal{S}^2_+}
\bullet\ R_{m^*,g}^{(\mu)}(\vect{r},\vect{\Omega})drd\Omega\equiv
\left(R_{m^*,g}^{(\mu)},\bullet\right)_{\mathcal{D}^2,\mathcal{S}^2_+}
\end{equation}
to yield radially-integrated cross sections for each layer $\mathcal{D}^1_j$
\begin{align}
\widetilde{\Sigma}_{t,g,j}^{(m^*,m)}&\equiv\sum_{i=1}^I \left(R_{m^*,g},R_{m,g}\right)_{\mathcal{D}^2_i,\mathcal{S}^2_+} \Sigma_{t,g,i,j}\,, \\
\widetilde{\Sigma}_{s,0,g'\rightarrow g,j}^{(m^*,m)}&\equiv\sum_{i=1}^I
\sum_{\ell=0}^L\!\sum_{\substack{k=-\ell\\\text{even }\ell+k}}^{\ell}
\frac{2\ell+1}{2\pi}\left(R_{m^*,g},Y_{\ell,k} \int_{\mathcal{S}^2_+}Y_{\ell,k}R_{m,g'}d\Omega'\right)_{\mathcal{D}^2_i,\mathcal{S}^2_+} \Sigma_{s,\ell,g'\rightarrow g,i,j}\,,\\
\widetilde{\Sigma}_{s,1,g'\rightarrow g,j}^{(m^*,m)}&\equiv\sum_{i=1}^I
\sum_{\ell=1}^L\!\sum_{\substack{k=-\ell\\\text{odd }\ell+k}}^{\ell}
\frac{2\ell+1}{2\pi}\left(R_{m^*,g},Y_{\ell,k} \int_{\mathcal{S}^2_+}Y_{\ell,k}R_{m,g'}d\Omega'\right)_{\mathcal{D}^2_i,\mathcal{S}^2_+} \Sigma_{s,\ell,g'\rightarrow g,i,j}\,,
\end{align}
which appear in the $z$-dependent cross sections
\begin{align}
\begin{split}
\widetilde{\Sigma}_{t,g}^{(m^*,m)}(z)&\equiv\sum_{j=1}^J \bm{1}_{\mathcal{D}^1_j}(z)\widetilde{\Sigma}_{t,g,j}^{(m^*,m)}+\left(
R_{m^*,g},\sqrt{1-\mu^2}\left(\cos(\omega)\frac{\partial}{\partial x}+\sin(\omega)\frac{\partial}{\partial y}
\right)R_{m,g}\right)_{\mathcal{D}^2,\mathcal{S}^2_+}\,,\end{split} \\
\widetilde{\Sigma}_{s,p,g'\rightarrow g}^{(m^*,m)}(z)&\equiv\sum_{j=1}^J \bm{1}_{\mathcal{D}^1_j}(z)\widetilde{\Sigma}_{s,p,g'\rightarrow g,j}^{(m^*,m)}\,.
\end{align}
Introducing yet another (coaxial) streaming coefficient for convenience,
\begin{equation}
s_{\oneD,g}^{(m^*,m)}\equiv\left(R_{m^*,g},\mu R_{m,g}\right)_{\mathcal{D}^2,\mathcal{S}^2_+}
\end{equation}
the axial operator for group $g$ becomes succinctly
\begin{equation}
\begin{split}
\mathrm{B}^{(m^*,m)}_{\oneD,g}Z_{m}&\equiv
s^{(m^*,m)}_{\oneD,g}\dir\frac{\partial}{\partial z}Z_{m,g}(z,\dir)
+\widetilde{\Sigma}_{t,g}^{(m^*,m)}(z)Z_{m,g}(z,\dir)
-\frac{1}{2}
\sum_{p=0}^{1} \dir^{p}
\sum_{g'=1}^G \widetilde{\Sigma}_{s,p,g'\rightarrow g}^{(m^*,m)}(z)
Z_{m,g}^p(z)
\,
\end{split}
\end{equation}
where $p$ and $Z^p_{m,g}$ are as in the previous section,
from which it follows that the reduced, 1D equation for the axial mode  $Z_{m^*,g}$ is
\begin{equation}
\label{eq:prog-1D}
\begin{split}
\mathrm{B}_{\oneD}^{(m^*,m^*)}Z_{m^*}&=\widetilde{q}_{\oneD}^{(m^*)}-\sum_{m\neq m^*}^M \mathrm{B}_{\oneD}^{(m^*,m)}Z_{m} \\
&=\widetilde{r}_{\oneD}^{(m^*)}
\end{split}
\end{equation}
in which
\begin{equation}
\widetilde{q}^{(m^*)}_{\oneD,g}(z,\dir)\equiv\sum_{f=1}^F \left(R^{(\mu)}_{m^*,g},q^{\twoD(\mu)}_{f,g}\right)_{\mathcal{D}^2,\mathcal{S}^2_+} q^{\oneD}_{f,g}(z,\dir)\,.
\end{equation}
As in previous sections, the streaming coefficient can be canceled by substituting
\begin{equation}
\begin{split}
\widetilde{\Sigma}_{t,g}^{(m^*,m)}(z)&\leftarrow\widetilde{\Sigma}_{t,g}^{(m^*,m)}(z)/s_{\oneD,g}^{(m^*,m^*)}\,, \\
\widetilde{\Sigma}_{s,p,g'\rightarrow g}^{(m^*,m)}(z)&\leftarrow\widetilde{\Sigma}_{s,p,g'\rightarrow g}^{(m^*,m)}(z)/ s_{\oneD,g}^{(m^*,m^*)}\,, \\
\widetilde{r}^{(m^*)}_{\oneD,g}(z)&\leftarrow\widetilde{r}^{(m^*)}_{\oneD,g}(z)/s_{\oneD,g}^{(m^*,m^*)}\,.
\end{split}
\end{equation}
Despite the fact that the axial mode $Z_{m^*}$ is independent of $\mu$, Equation \ref{eq:prog-1D} can be viewed as an $S_N$ (discrete ordinates) discretization of 1D neutron transport with a degenerate quadrature consisting of two ordinates of equal weight, at angles $\{-1,+1\}$.
Alternatively, by recognizing that
the even and odd axial fluxes are proportional to the 1D scalar flux and current,
Equation \ref{eq:prog-1D} can be rearranged into a form similar to the well-known $P_1$ equations; this possibility is discussed in \ref{app:1D-P1}.

\subsection{Radial or Axial Energy Dependence}
\label{sec:radial-axial-energy}
In the above discussion, we granted that the energy dependence would be assigned to both the radial and axial modes. In the multigroup discretization, this amounts to a group-wise 2D/1D decomposition. This is viable, and in fact has been the exclusive purview of previous work on spatial separation of neutron diffusion \cite{Senecal2019, Alberti2020, Prince2020b}. However, as separation of variables is a defining feature of PGD, it is clear that the energy dependence could alternatively be assigned to either the radial or axial modes and separated out of the other. Especially as the number of energy groups increases, this may reduce the computational cost of the ROM substantially. 
Moreover, as computing the radial modes will generally be much more burdensome than the same for the axial modes, the savings may be even greater in the case where the radial equations are rendered energy-independent. 

In any case, deriving the energy independent equations from their multigroup counterparts is straightforward, especially as the former can be considered a degenerate case of the latter where $G=1$ (that is, there is only one group). The coefficients for this single group, denoted $\mathcal{G}$ for clarity, can be found by integrating each energy dependent term over $[E_G,\,E_0]$. 
This is, in the multigroup approximation,
equivalent to a summation over $g=1\ldots G$, with two caveats: first, group fluxes are conventionally defined as integral terms, so that
\begin{equation}
\psi_g(\vec{r},\vec{\Omega})=\int_{E_{g}}^{E_{g-1}}\psi(\vec{r},\vec{\Omega},E) dE\,,
\end{equation}
and so there is, strictly speaking,  
no means by which to evaluate inner products of form $\int_{E_{g}}^{E_{g-1}}\psi\psi^* dE$ given a pair of group fluxes, $\psi_g$ and $\psi_g^*$. This is trivially overcome by assuming the flux is constant within a group, such that $\psi_g=\bar{\psi}_g\Delta_g$ (likewise $\psi^*_g$)
where $\bar{\psi}_g$ is the average flux and $\Delta_g=E_{g-1}-E_{g}$ is the group width, allowing the previous inner product to be written  $\Delta_g^{-1}\psi_g\psi_g^*$.
The second nuance is that the energy scales involved in neutron slowing down are inherently logarithmic, as the physics are that of a particle scattering many times and losing a fraction of its energy each time. Accordingly, a more appropriate inner product space could be given in lethargy $u$ rather than energy $E$, defined as $u=\ln(E_0/E)$. 
Enacting this change is accomplished by simply redefining the group width in units of lethargy, $\Delta_g\leftarrow \ln(E_{g-1}/E_{g})$.

Given these preliminaries, we can succinctly express the one-group coefficients as
\begin{align}
s_{\mathrm{ND},\mathcal{G}}^{(m^*,m)}&\equiv\sum_{g=1}^G \Delta_g^{-1} s^{(m^*,m)}_{\mathrm{ND},g}\,, \\
\widetilde{\Sigma}_{t,\mathcal{G}}^{(m^*,m)}&\equiv\sum_{g=1}^G \Delta_g^{-1}\widetilde{\Sigma}_{t,g}^{(m^*,m)}\,, \\
\widetilde{\Sigma}_{s,\ell,k,\mathcal{G}}^{(m^*,m)}&\equiv\sum_{g=1}^G \Delta_g^{-1}\sum_{g'=1}^G \widetilde{\Sigma}_{s,\ell,k,g'\rightarrow g}^{(m^*,m)}\,,
\end{align}
and likewise the one-group source,
\begin{equation}
q^{(m^*)}_{\mathrm{ND},\mathcal{G}}\equiv\sum_{g=1}^G \Delta_g^{-1} q^{(m^*)}_{\mathrm{ND},g}\,,
\end{equation}
where the subscript $\mathrm{ND}$ refers to the sub-model from which the energy variable is being separated, either $\twoD$, $\oneD(\mu)$, $\twoD(\mu)$, or $\oneD$ as in Sections \ref{sec:theory-axial-polar-pgd} and \ref{sec:theory-axial-pgd}. Beyond this, the corresponding mode ($R_m$, $Z^{(\mu)}_m$, $R^{(\mu)}_m$, or $Z_m$) is redefined as independent of $g$, as are the source components $q^{\mathrm{ND}}_f$.
With these modifications, the multigroup models---Equations \ref{eq:prog-2D}, \ref{eq:prog-1D-mu}, \ref{eq:prog-1D-mu}, and \ref{eq:prog-2D-mu}---simplify to one-group, but otherwise retain all aforementioned particularities, such as the absence of $\mu$ or polar anisotropy of the total cross section. 

\subsection{Proper Generalized Decomposition Algorithm}
\label{sec:alg-pgd-2D1D}
Given these 2D and 1D sub-models, where either or both (but not neither) is multigroup, the overall iterative procedure of PGD can be enacted simply.
To review, Progressive PGD begins with zero modes, $M=0$, and employs a greedy algorithm to find the first pair. To do, some nonlinear iteration is required between the axial and radial sub-models for $m^*=M=1$; here, as in most PGD literature, fixed-point (sometimes termed ``Alternating Directions'') iteration is chosen for its simplicity and demonstrated efficacy. Nonetheless, several enhancements to this basic procedure exist, see \cite{Senecal2017}, which may be advantageous in practice. Upon nonlinear convergence, this process then repeats, ``enriching'' the reduced basis---as in, incrementing $M$---until some numerical tolerance or other termination criterion is met. This procedure is detailed in Algorithm \ref{alg:progressive}, while the specific iterative solvers and settings are given in Section \ref{sec:sundries}.

\begin{algorithm}
	\caption{Axial/-polar Proper Generalized Decomposition}
	\SetKw{KwElse}{else}
	\DontPrintSemicolon
	\label{alg:progressive}
	\tcp{$R_{m^*}$ and $Z_{m^*}$ denote radial and axial modes; superscript $(\mu)$ is omitted. Either}
	\tcp{or both may be energy-dependent (multigroup).}
	$\psi\gets0$ \\
	\For(\tcp*[f]{enrichment iteration}){$m^*=1\ldots M$}{
		$R_{m^*}\gets 1$ \\
		\While(\tcp*[f]{nonlinear iteration}){\texttt{residual} $>$ \texttt{tolerance}}{
			\tcp{Normalize the radial mode for numerical stability}
			\If{$R_{m^*}$ and $Z_{m^*}$ are both multigroup}{
				\lFor{$g=1\ldots G$}{$R_{m^*,g}\gets R_{m^*,g}/\lVert R_{m^*,g}\rVert_{L^2}$}
			}
			\Else{
				$R_{m^*}\gets R_{m^*}/\lVert R_{m^*}\rVert_{L^2}$ \\
			}
			\tcp{Solve the relevant submodels for axial-polar or axial PGD}
			$R_{m^*}\gets$ solution of Equation \ref{eq:prog-2D} or \ref{eq:prog-2D-mu} \\
			$\texttt{residual},\ Z_{m^*}\gets$ initial residual and solution of Equation \ref{eq:prog-1D-mu} or \ref{eq:prog-1D} \\
		}
		$\psi\gets\psi+R_{m^*}\otimes Z_{m^*}$
	}
\end{algorithm}

\subsection{Similarity to Conventional 2D/1D Methods}
\label{sec:pgd-vs-2D1D}
Despite their independent origin, these 2D/1D PGD models share some similarities to existing 2D/1D methods (as discussed briefly in Section \ref{sec:prev-work}). Namely, taking the Galerkin projection with respect to the transverse spatial dimension(s)---that is, applying Equations \ref{eq:proj-2D}, \ref{eq:proj-1D-mu}, \ref{eq:proj-2D-mu}, or \ref{eq:proj-1D}---is essentially equivalent to homogenization by transverse integration as practiced in many 2D/1D methods, except that we introduce a (non-uniform) test function \cite{Collins2016}. Likewise, anisotropic leakage terms in the total cross section also appear in MPACT's 2D/1D method, through a procedure termed ``transverse leakage splitting'' \cite{Collins2016, Jarrett2018}.
Finally, although we do not solve it as such here, the 1D equation of axial PGD is equivalent to a diffusion or $P_1$ model (albeit one that is not limited to $L\leq1$ scattering, see \ref{app:1D-P1}), which is the conventional axial discretization in many 2D/1D methods \cite{Kelley2015}.

The main differences, meanwhile, are first that we do not approximate the angular dependence of the flux in either the transverse leakage terms or the other inner products that appear in the cross sections. This is commonly done in conventional 2D/1D methods, where the transverse leakage is assumed to be isotropic (or approximated as a truncated series) and the cross section ``homogenization'' is applied using the scalar rather than angular flux (which Jarrett \cite{Jarrett2018} and Faure et al. \cite{Faure2018} identify as an appreciable source of error). Secondly, PGD seeks multiple modes, which allows it to converge to the 3D solution (hypothetically, though in practice PGD's convergence may be prohibitively slow), while other 2D/1D methods seek only one---in effect, approximating the solution as rank-one separable. Thirdly, and relatedly, conventional 2D/1D methods subdivide the domain into several planes and unit-cells (commonly pin-cells) in order to localize the 2D/1D decomposition. This could be done in PGD as well, but is not strictly necessary since (supposing that PGD converges reasonably quickly) the number of modes can simply be increased until the solution is adequately resolved throughout the domain, as measured by some error indicator. Fourthly, while 2D/1D methods have been demonstrated with both 1D transport and diffusion (analogous to our axial-polar and axial submodels), 
the angular dependence on $\mu$ is shared between the 2D and 1D models in the former case. 
This could also be considered in PGD (analogously to the case where both axial and radial modes are multigroup), but we find it more natural to relegate $\mu$ to either the 1D or 2D domain, but not both, in the axial-polar PGD. Finally, we here use the same discretization in the 2D and 1D sub-models---such that we can easily reconstruct the corresponding 3D solution---whereas many 2D/1D models employ the Method of Characteristics in the 2D planes and nodal diffusion in the 1D axes. That said, this is not an essential feature of either method: one could conceivably apply different discretizations to the 2D and 1D sub-models of PGD just as easily as one could use the same discretization for both in traditional 2D/1D methods (see Jarrett's ``1D/1D $S_N$'' method \cite{Jarrett2018}, for example).

\section{Numerical Results}
\label{sec:results-2D1D}
To characterize performance, we begin by solving the first Takeda benchmark \cite{Takeda1991}, a small
homogenized LWR core, and an infinite lattice of UO\textsubscript{2} pin-cells, as specified in the (modified) C5G7 benchmark \cite{Lewis2005}. 
Performance is primarily measured by the $L^2$ norm of the error of the angular and scalar fluxes, as compared to a full-order reference solution. Moreover, the ideal decomposition in the $L^2$ norm 
given by the SVD of the reference solution
is shown for comparison.\footnote{More specifically, a separate reference decomposition is computed for each of the six candidate ROMs, either by the SVD of the 3D multigroup angular flux or, in the case where both axial and radial modes are multigroup, by the SVD of the 3D angular flux in each individual group. This yields several decompositions of the angular flux, which are subsequently integrated in angle to produce corresponding decompositions of the scalar flux.
}
Results are given both as a function of rank $M$ and, for PGD, runtime---that is, the real or ``wall-clock'' time elapsed since computation began, not including any requisite pre- and post-processing. The same is not provided for the full-order (3D) calculation, as performance-critical optimizations---most notably, diffusion preconditioning and parallelism---have yet to be implemented in either the full- or reduced-order models, rendering any quantitative comparison of the two premature and perhaps misleading. (That said, the potential for large computational savings is evident by the definition of the ROM and was qualitatively observed in practice.) 
However, as each of the six PGD ROMs employs the same iterative procedures, numerical tolerances, and overall algorithm, comparison among them appears apt, given any future enhancements---such as PGD's ``update'' step \cite{Nouy2010,Dominesey2022}, Minimax projection \cite{Nouy2010,Dominesey2022b,Dominesey2022}, or a Newton or quasi-Newton nonlinear solver \cite{Ammar2006,Dennis1977}---stand to benefit each. Accordingly, while we do not expect the reported runtimes to be representative in absolute terms (except as an upper bound to be improved upon), the performance of each relative to the others may hold (at least roughly) fixed even as the algorithm and underlying implementation are improved. As such, these results are practically useful, especially in informing the direction of future research.

\subsection{Numerical Sundries}
\label{sec:sundries}
\begin{figure}[htbp]
	\centering
	\begin{subfigure}{0.32\textwidth}
		\begin{tikzpicture}[scale=1.8,tdplot_main_coords]
		
		
		\def\rvec{1}
		\def\thetavec{30}
		\def\phivec{60}
		
		\def\nAzi{16}
		\def\nPol{4}
		\pgfmathsetmacro{\wedge}{360/\nAzi}
		\pgfmathsetmacro{\from}{\wedge/2}
		
		\coordinate (O) at (0,0,0);
		\draw[thick,->] (0,0,0) -- (1.2,0,0) node[anchor=south east]{$x'$};
		\draw[thick,->] (0,0,0) -- (0,1.2,0) node[anchor=south west]{$y'$};
		\draw[thick,<->] (0,0,-1.3) -- (0,0,1.2) node[anchor=south]{$z'$};
		
		\def\angEl{30} 
		\begin{scope}[tdplot_screen_coords]
		\DrawLongitudeCircleBoth[1]{-20}{gray,dashed,fill=black,fill opacity=0.15}
		\end{scope}	
		
		\begin{scope}[tdplot_screen_coords]
		\DrawLongitudeCircleBoth[1]{-110}{gray,dashed,fill=black,fill opacity=0.0}
		\end{scope}
		
		\tdplotdrawarc[draw=none,gray,dashed,thick,fill=black,fill opacity=0.15]{(O)}{1}{0}{90}{}{};
		\fill[black,opacity=0.15] (1,0,0) -- (0,0,0) -- (0,1,0) -- cycle;
		
		\def\wgts{{0.3626837833783620,0.3137066458778873,0.2223810344533745,0.1012285362903763}}
		\def\mus{{0.1834346424956498,0.5255324099163290,0.7966664774136267,0.9602898564975363}}

		\foreach \s in {1,-1}
		\foreach \j in {\nPol,...,1} {
			\pgfmathsetmacro{\mu}{\s*\mus[\j-1]}
			\pgfmathsetmacro{\wgt}{\wgts[\j-1]}
			\pgfmathsetmacro{\proj}{sqrt(1-\mu*\mu)}
			\tdplotdrawarc[gray,thick]{(0,0,\mu)}{{\proj}}{0}{90}{}{};
			\foreach \i in {1,...,4} {
				\pgfmathsetmacro{\omega}{\i*\wedge-\from}
				\ifnum \j=\nPol
				\ifnum \s=1
				\begin{scope}[tdplot_screen_coords]
				\DrawLongitudeCircleBoth[1]{\omega-110}{gray,thick}
				\end{scope}
				\else\fi
				\else\fi
				\shade[ball color=black] ({\proj*cos(\omega)},{\proj*sin(\omega)},\mu) circle ({\wgt*4.5pt});
			}
		}
		\end{tikzpicture}
		\caption{}
	\end{subfigure}
	\hfill
	\begin{subfigure}{0.32\textwidth}
		\begin{tikzpicture}[scale=1.8,tdplot_main_coords]
		
		
		\node at (0,0,-1.3) {};
		
		\def\rvec{1}
		\def\thetavec{30}
		\def\phivec{60}
		
		\def\nAzi{16}
		\def\nPol{4}
		\pgfmathsetmacro{\wedge}{360/\nAzi}
		\pgfmathsetmacro{\from}{\wedge/2}
		
		\coordinate (O) at (0,0,0);
		\draw[thick,->] (0,0,0) -- (1.2,0,0) node[anchor=south east]{$x'$};
		\draw[thick,->] (0,0,0) -- (0,1.2,0) node[anchor=south west]{$y'$};
		\draw[thick,<->] (0,0,-1.2) -- (0,0,1.2) node[anchor=south]{$z'$};
		
		\def\angEl{30} 
		\begin{scope}[tdplot_screen_coords]
		\DrawLongitudeCircleBoth[1]{-20}{dashed,draw=gray,fill=red,fill opacity=0.15}
		\end{scope}
		
		\def\wgts{{0.3626837833783620,0.3137066458778873,0.2223810344533745,0.1012285362903763}}
		\def\mus{{0.1834346424956498,0.5255324099163290,0.7966664774136267,0.9602898564975363}}
		
		\foreach \s in {-1,1} {
			\foreach \j in {1,...,\nPol} {
				\pgfmathsetmacro{\mu}{\s*\mus[\j-1]}
				\pgfmathsetmacro{\wgt}{\wgts[\j-1]}
				\pgfmathsetmacro{\proj}{sqrt(1-\mu*\mu)}
				\draw[red,thick] (0,0,\mu) -- (0,\proj,\mu);
				\shade[ball color=red] (0,0,\mu) circle ({\wgt*4.5pt});
			}
		}
		\tdplotdrawarc[gray,dashed,draw=gray,fill=blue,fill opacity=0.15]{(O)}{1}{0}{90}{}{};
		\fill[blue,opacity=0.15] (1,0,0) -- (0,0,0) -- (0,1,0) -- cycle;
		\foreach \i in {1,...,4} {
			\pgfmathsetmacro{\omega}{\i*\wedge-\from}
			\draw[blue,thick] (O) -- ({cos(\omega)},{sin(\omega)},0);
			\shade[ball color=blue] ({cos(\omega)},{sin(\omega)},0) circle (1.5pt);
		}
	
		\end{tikzpicture}
		\caption{}
	\end{subfigure}
	\hfill
	\begin{subfigure}{0.32\textwidth}
		\begin{tikzpicture}[scale=1.8,tdplot_main_coords]
		
		
		\node at (0,0,-1.3) { };
		
		\def\rvec{1}
		\def\thetavec{30}
		\def\phivec{60}
		
		\def\nAzi{16}
		\def\nPol{4}
		\pgfmathsetmacro{\wedge}{360/\nAzi}
		\pgfmathsetmacro{\from}{\wedge/2}
		
		\coordinate (O) at (0,0,0);
		\draw[thick,red] (0,0,-1) -- (0,0,1);
		\draw[thick,->] (0,0,0) -- (1.2,0,0) node[anchor=south east]{$x'$};
		\draw[thick,->] (0,0,0) -- (0,1.2,0) node[anchor=south west]{$y'$};
		\draw[thick,->] (0,0,1) -- (0,0,1.2) node[anchor=south]{$z'$};
		\draw[thick,->] (0,0,-1) -- (0,0,-1.2);
		
		\def\angEl{30} 
		\begin{scope}[tdplot_screen_coords]
		\DrawLongitudeCircle[1]{-20}{gray,dashed,fill=blue,fill opacity=0.15}
		\end{scope}	
		\fill[blue,opacity=0.15] (0,0,1) -- (0,0,0) -- (0,1,0) -- cycle;
		
		\begin{scope}[tdplot_screen_coords]
		\DrawLongitudeCircle[1]{-110}{gray,dashed,fill=black,fill opacity=0.0}
		\end{scope}
		
		\tdplotdrawarc[draw=gray,dashed,fill=blue,fill opacity=0.15]{(O)}{1}{0}{90}{}{};
		\fill[blue,opacity=0.15] (1,0,0) -- (0,0,0) -- (0,1,0) -- cycle;
		
		\def\wgts{{0.3626837833783620,0.3137066458778873,0.2223810344533745,0.1012285362903763}}
		\def\mus{{0.1834346424956498,0.5255324099163290,0.7966664774136267,0.9602898564975363}}
		
		\foreach \j in {1,...,\nPol} {
			\pgfmathsetmacro{\mu}{\mus[\j-1]}
			\pgfmathsetmacro{\wgt}{\wgts[\j-1]}
			\pgfmathsetmacro{\proj}{sqrt(1-\mu*\mu)}
			\tdplotdrawarc[blue,thick]{(0,0,\mu)}{{\proj}}{0}{90}{}{};
			\foreach \i in {1,...,4} {
				\pgfmathsetmacro{\omega}{\i*\wedge-\from}
				\ifnum \j=1
				\begin{scope}[tdplot_screen_coords]
				\DrawLongitudeCircle[1]{\omega-110}{blue,thick}
				\end{scope}
				\else\fi
				\shade[ball color=blue] ({\proj*cos(\omega)},{\proj*sin(\omega)},\mu) circle ({\wgt*4.5pt});
			}
		}
		
		\shade[ball color=red] (0,0,-1) circle (1.5pt);
		\shade[ball color=red] (0,0,+1) circle (1.5pt);
		
		\end{tikzpicture}
		\caption{}
	\end{subfigure}
	\caption{Angular quadratures of full- (a) and reduced-order 
		(c/b) models by axial/-polar 
		Proper Generalized Decomposition.}
	\label{fig:quadratures}
\end{figure}
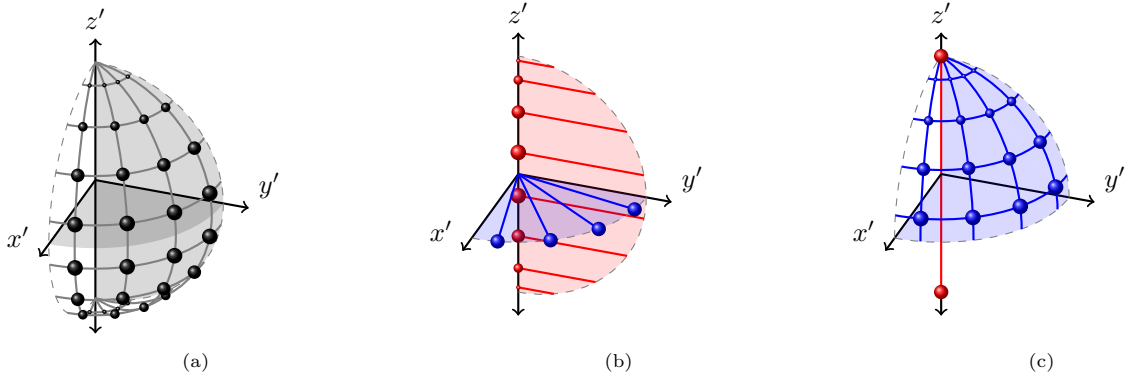

Each numerical experiment is performed with an $S^4_4$ product Gauss-Legendre-Chebyshev angular quadrature \cite{AbuShumays1977}---that is, 4 polar and 4 azimuthal angles per octant, 128 in total. As visualized in Figure \ref{fig:quadratures}, this quadrature---like any
defined as the tensor product of polar and azimuthal quadratures---allows the angular flux to be reconstructed from the axial-polar PGD.
Likewise, finite elements formed by (spatial) tensor products, such as the linear discontinuous Lagrangian elements employed herein, allow the 3D flux to be reconstructed straightforwardly from 2D and 1D modes.
Further, because each benchmark is specified as a criticality (that is, $k$-eigenvalue) problem rather than fixed-source problem, a uniform source is assumed throughout the core; this is typically a crude assumption of the true fission source, but a reasonable guess by which to initialize some iterative eigenvalue procedure.
In this sense, we consider the present case sufficiently representative of a simple use-case for 2D/1D methods.

Each reference solution is computed by the Generalized Minimal Residual (GMRES) method with a tolerance (on the $\ell^2$ norm of the residual) of $10^{-6}$ times the $\ell^2$ norm of the uncollided flux, restarted after Krylov subspace reaches a maximum size of 5---usually denoted as GMRES(5). This iteration is preconditioned by block Gauss-Seidel in energy, wherein each diagonal (within-group) block is solved by GMRES(15), terminated when the residual falls below either $10^{-6}$ or $10^{-2}$ times its initial value (hereafter referred to as absolute and relative tolerances) or after a maximum of 250 iterations.

The PGD ROMs, meanwhile, are run for a fixed number of 30 modes. This is for testing purposes, as in most other scenarios one would select an error indicator and tolerance rather than $M$ itself. 
The alternating, nonlinear iteration is terminated after 10 iterations or when the normalized, initial, 1D or 1D($\mu$) residual falls below $10^{-2}$. 
More specifically, the latter metric refers to,  in the 1D case, the $\ell^2$ norm of residual of the system
$\mathrm{B}^{(m^*,m^*)}_{\oneD}Z_{m^*}=\widetilde{r}^{(m^*)}_{\oneD}$ (before it is solved) divided by that of $\widetilde{r}^{(m)}_{\oneD}$---or more concisely, $\lVert r^{(m^*+1)}_{\oneD}\rVert_{\ell^2}/\lVert r^{(m^*)}_{\oneD}\rVert_{\ell^2}$. The 1D($\mu$) indicator is identical, though based on Equation \ref{eq:prog-1D-mu} instead of \ref{eq:prog-1D}. Each submodel is solved by GMRES(25) with a maximum of 100 iterations and absolute and relative tolerances of $10^{-5}$ and $10^{-3}$. As in the full-order model, this is preconditioned by block Gauss-Seidel in energy, with each diagonal block being solved by GMRES(25) with a maximum of 500 iterations absolute and relative tolerances of $10^{-4}$ and $10^{-6}$. This applies to both multigroup and energy-independent submodels, as the latter are implemented as degenerate (one-group) cases of the former.
The implementation itself is written in object-oriented C++ using version 9.1 of the \verb|deal.II| finite element library \cite{Bangerth2007, Bangerth2019}.

\subsection{Takeda Light Water Reactor}
\label{sec:takeda-lwr}
The first of the Takeda benchmarks \cite{Takeda1991}---a small LWR with octant symmetry representing the Kyoto University Critical Assembly (KUCA)---is depicted in Figure \ref{fig:extruded-takeda}, where the red (A1), blue (B1, B2, A2), and white volumes (C1, C2) correspond to the homogenized core, water reflector, and control rod channel respectively; this channel is filled with either a neutron absorber or void in the first (rodded) and second (unrodded) cases respectively. As per benchmark specifications, energy and spatial discretizations are given by two-group cross sections and a Cartesian grid of uniform 1 cm spacing.

\begin{figure}[htb]
	\centering
	\begin{subfigure}[t]{0.49\textwidth}
		\centering
		\includegraphics[width=0.7\textwidth]{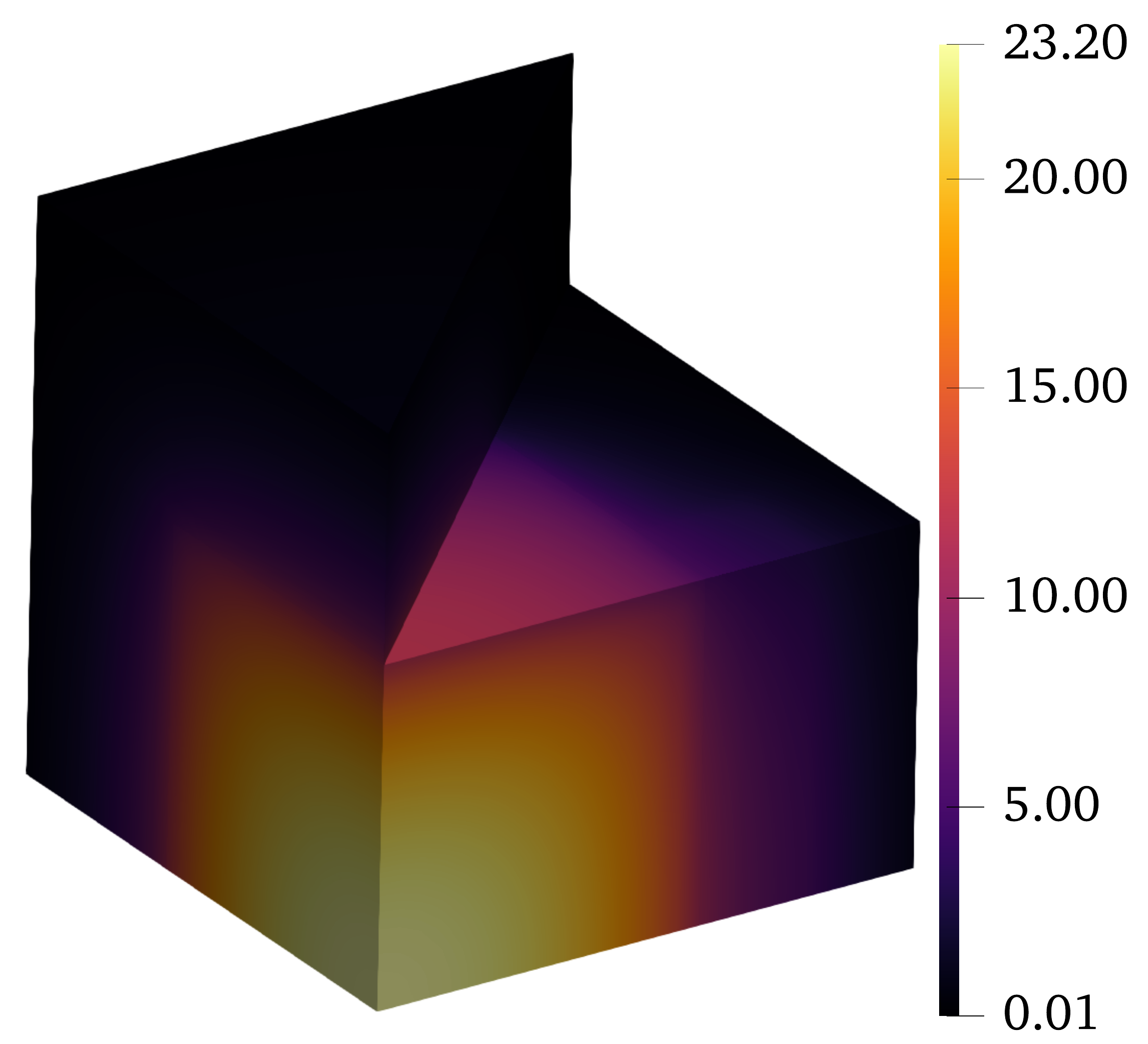}
		\caption{}
	\end{subfigure}
	\begin{subfigure}[t]{0.49\textwidth}
		\centering
		\includegraphics[width=0.7\textwidth]{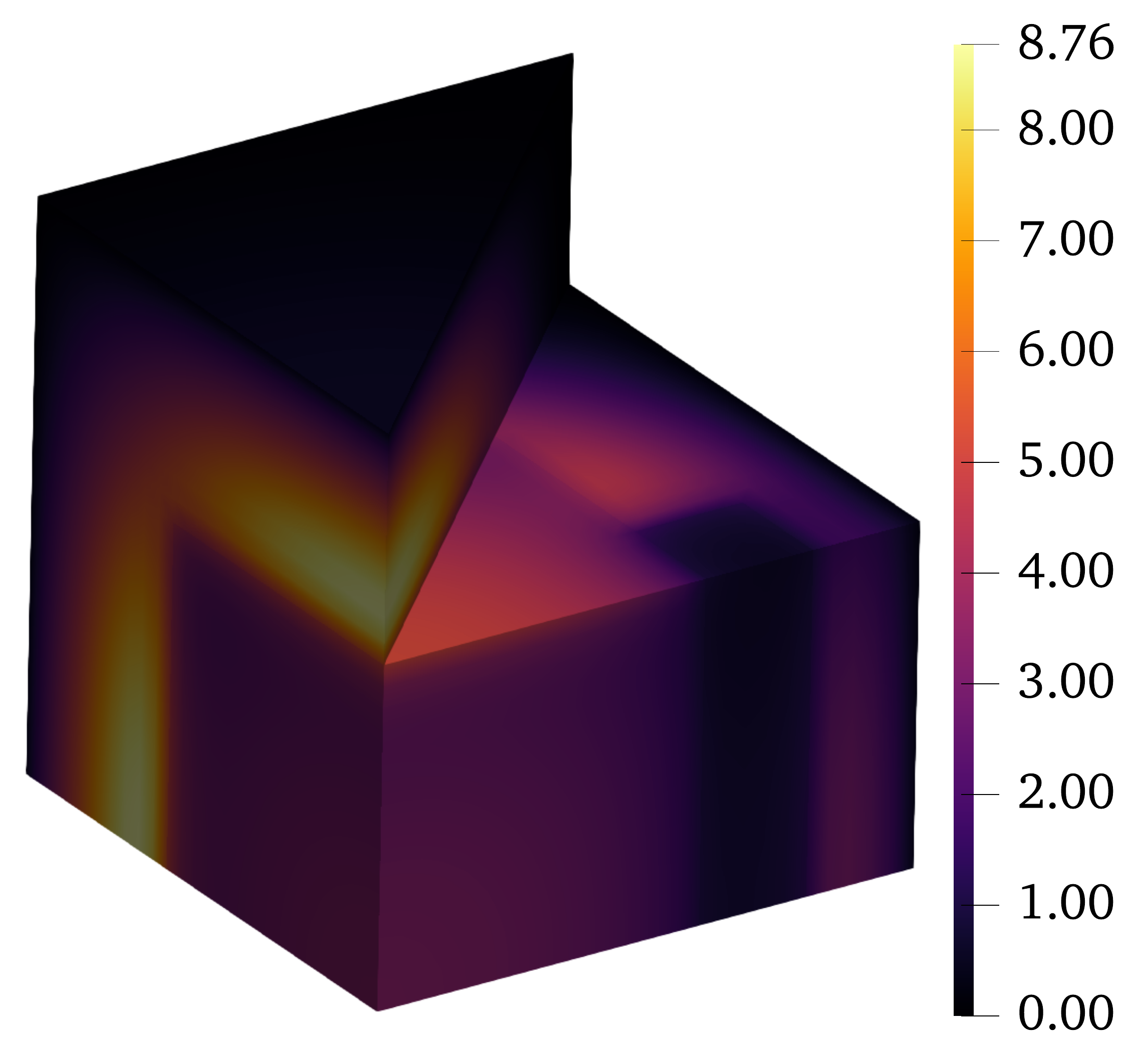}
		\caption{}
	\end{subfigure}
	\\
	\begin{subfigure}[t]{0.49\textwidth}
		\centering
		\includegraphics[width=0.7\textwidth]{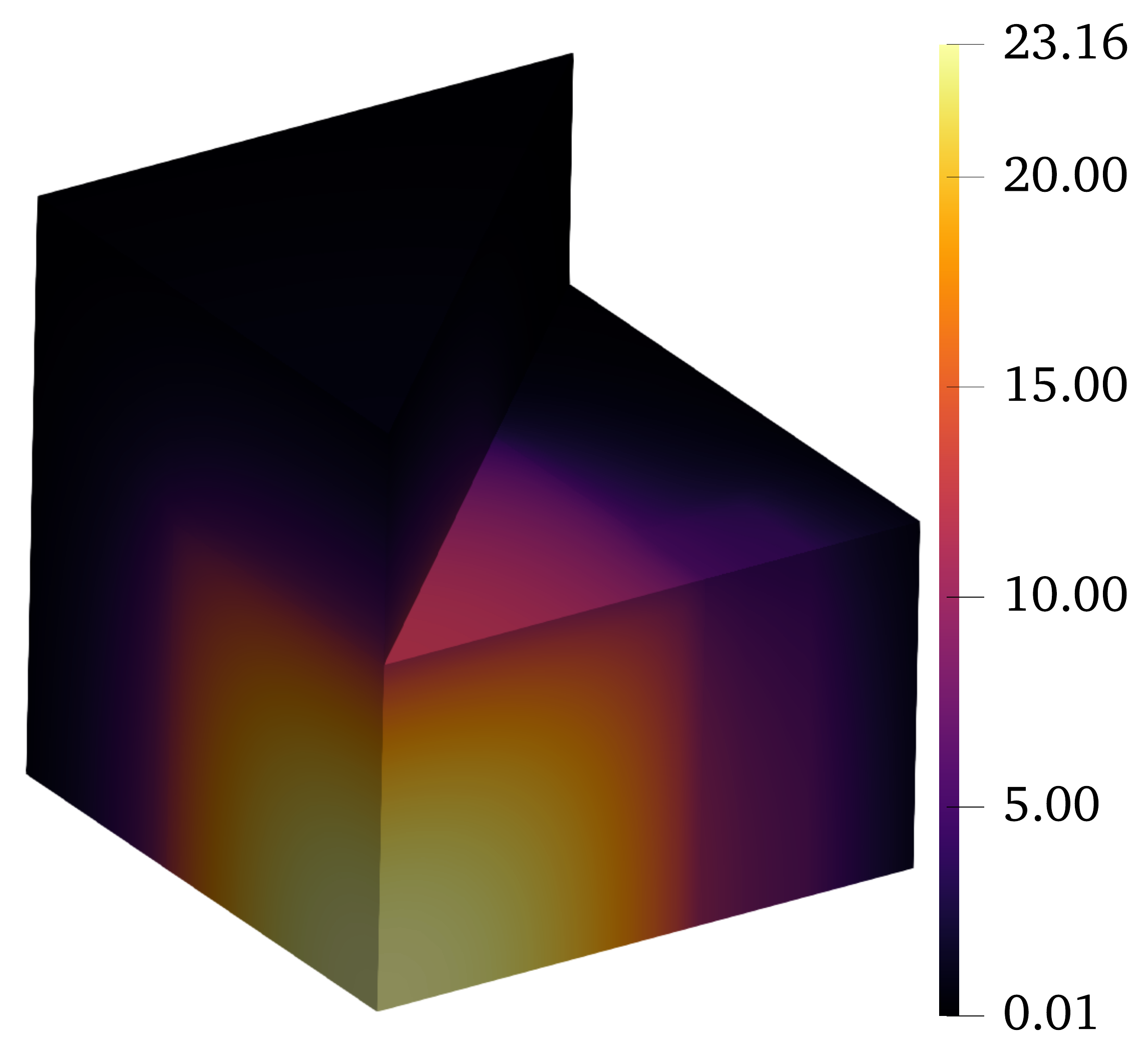}
		\caption{}
	\end{subfigure}
	\begin{subfigure}[t]{0.49\textwidth}
		\centering
		\includegraphics[width=0.7\textwidth]{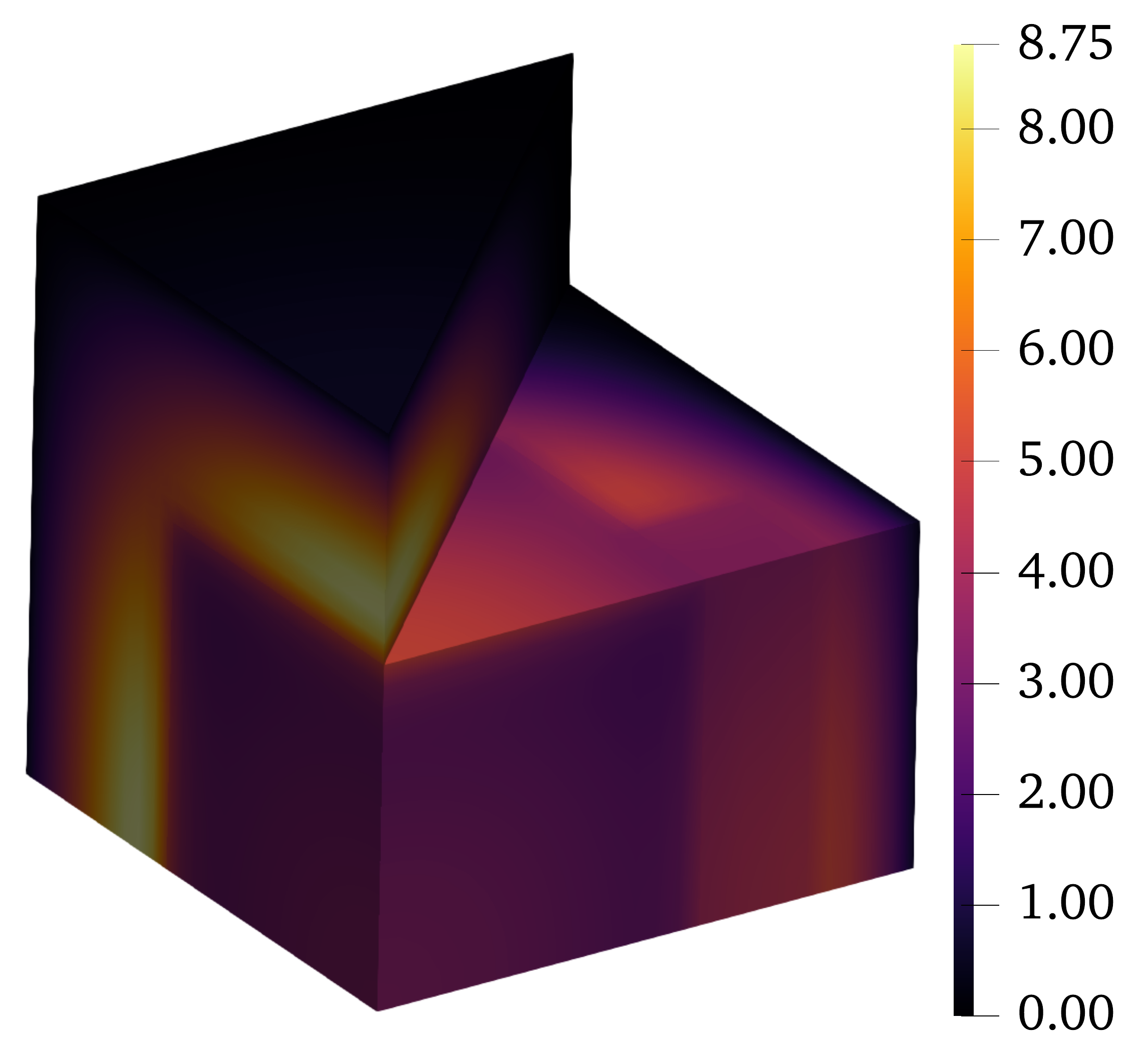}
		\caption{}
	\end{subfigure}
	\caption{Reference scalar fluxes for the 3D rodded (a,b) and unrodded (c,d) Takeda Light Water Reactor, in fast (a,c) and thermal (b,d) energy groups, $g=1,\,2$. The diagonal cutaway reveals the radial distribution near the rod channel.}
	\label{fig:takeda-3D}
\end{figure}

\begin{figure}
	\centering
	\includegraphics[width=0.915\textwidth]{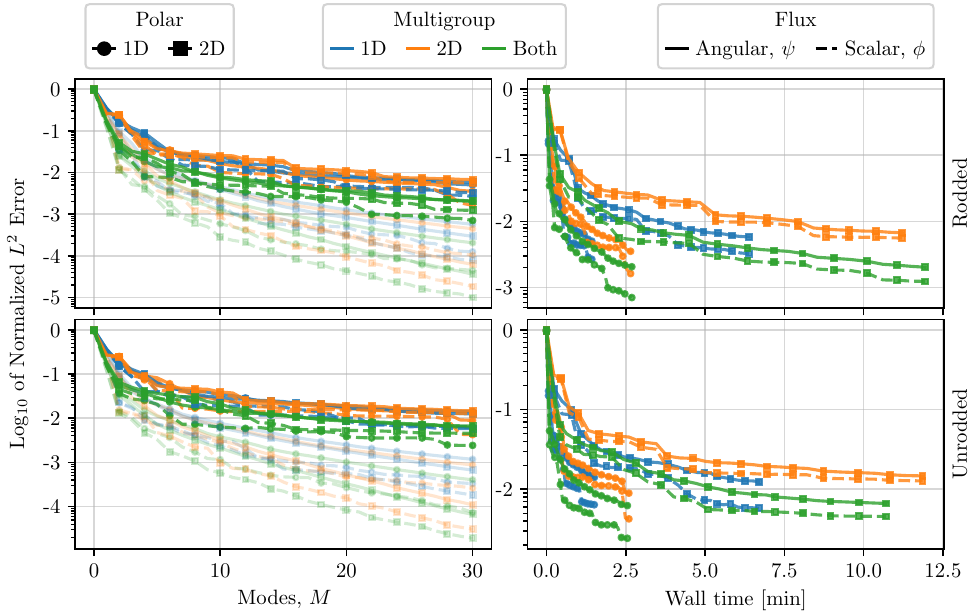}
	\caption{Convergence of 2D/1D Proper Generalized Decomposition (opaque lines) and Singular Value Decomposition (faint lines) for the Takeda Light Water Reactor.}
	\label{fig:2D1D-lwr}
\end{figure}

Given these parameters, the convergence of both SVD and PGD for each of the six decompositions in Figure \ref{fig:2D1D-lwr}. From this, we first observe that a low-rank approximation is possible and that PGD does indeed converge.
More specifically, by thirty modes the optimal decomposition of SVD is in all scenarios seen to achieve $L^2$ errors of the angular flux below $2\times10^{-3}$.
Moreover, we find the PGD, while not as rapidly convergent as the SVD, does yield an angular flux with an $L^2$ error below $2\%$ in each case at $M=30$. These results are somewhat superior in the group-wise PGD, decreasing to $0.1\%$ and $0.7\%$ for the rodded and unrodded configurations respectively.

However, while the convergence of each PGD scheme with increasing $M$ is comparable, the same does not apply with respect to runtime. Namely, the 2D($\mu$)/1D PGDs---those which assign the polar angle to the 2D sub-model---take substantially longer to compute. This aligns with the expectation that the computational cost of 2D/1D PGD will often be dominated by that of calculating the 2D modes. As such, removing $\mu$ from this computation---in this case reducing the number of angles in the 2D plane by a factor of four---yields an appreciable speedup,
yet appears to have little adverse effect on the quality of the decomposition. 
Meanwhile, a second, less severe, discrepancy is observed between the three treatments of energy. As to be expected, for both 2D($\mu$)/1D and 2D/1D($\mu$) schemes, the fastest-running model (per mode) is that where only the 1D modes are multigroup, then that for 2D modes, and finally the group-wise case, where both are multigroup.
However, unlike the previous case, this additional runtime is offset by a commensurate increase in accuracy, such that the group-wise PGD generally achieves the lowest error for a given amount of runtime.
That said, as noted in Section \ref{sec:radial-axial-energy}, group-wise PGD may be most advantageous when there are few energy groups with markedly different features arising in each (as clearly observed in Figure \ref{fig:takeda-3D});
given this, the net benefit demonstrated here may be diminished or reversed in finer-group problems.
In any case, these numerical results confirm each of the six 2D/1D ROMs to be feasible, convergent, and effective for this LWR benchmark.

\subsection{C5G7 Pin-Cell}
\begin{figure}[h!tb]
	\centering
	\begin{tikzpicture}
	\newcommand{\scale}{5}
	\newcommand{\pitch}{1.26cm*\scale}
	\newcommand{\radius}{0.54cm*\scale}
	\newcommand{\key}{0.6cm}
	\filldraw[fill=blue] (\pitch/2+\key,\pitch/2) rectangle node[anchor=west] {\hspace{0.5cm}Water (moderator)} +(\key,-\key);
	\filldraw[fill=red] (\pitch/2+\key,\pitch/2-1.5*\key) rectangle node[anchor=west] {\hspace{0.5cm}Homogenized fuel/clad} +(\key,-\key);
	\node[inner sep=0pt] at 
	(\pitch/4,\pitch/4)
	{\includegraphics[width=\pitch/2]{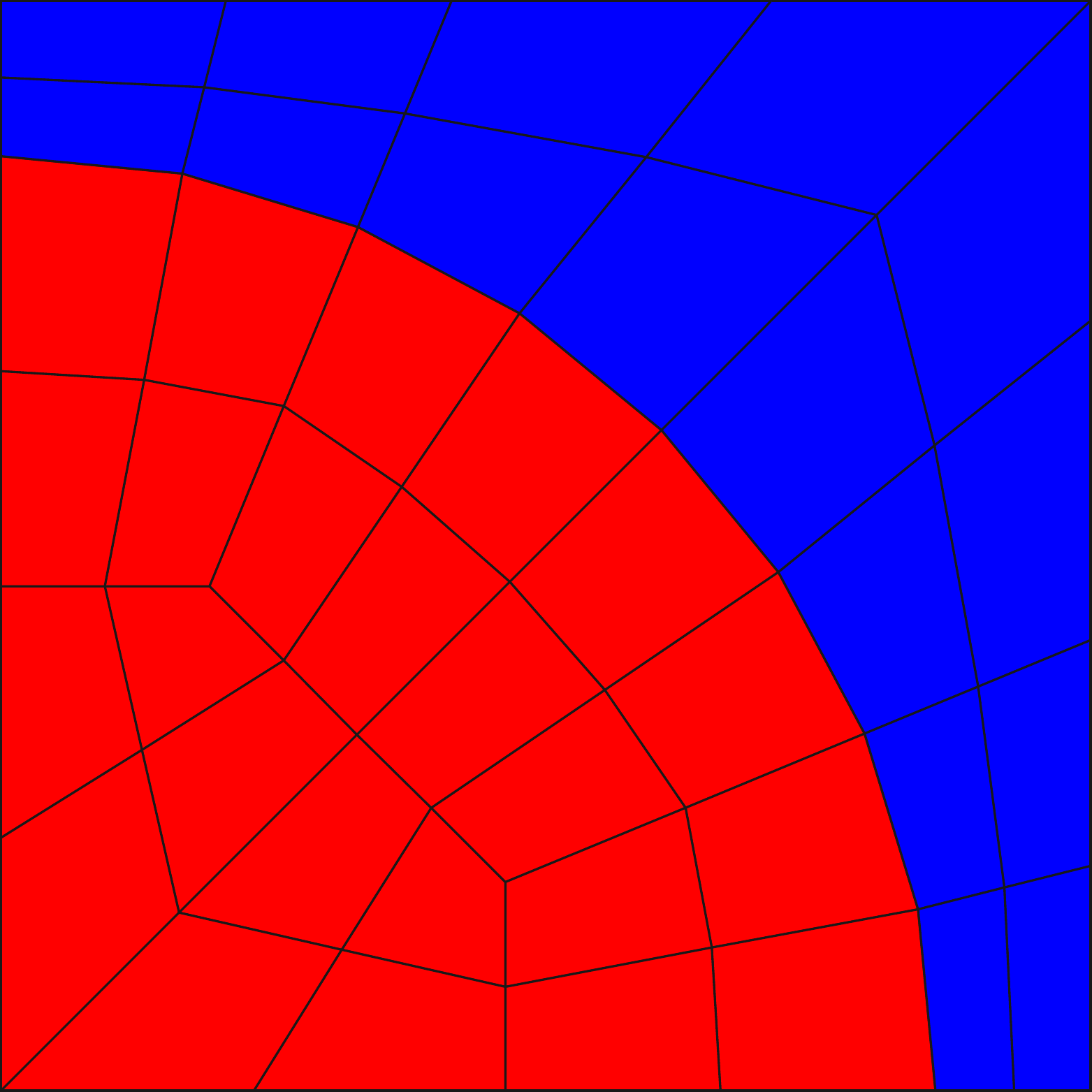}};
	\end{tikzpicture}
	\caption{Radial mesh of C5G7 UO\textsubscript{2} pin-cell, 5:1 scale.}
	\label{fig:pin-2D}
\end{figure}

\begin{figure}
	\centering
	\begin{subfigure}[t]{0.32\textwidth}
		\centering
		\includegraphics[width=0.425\textwidth]{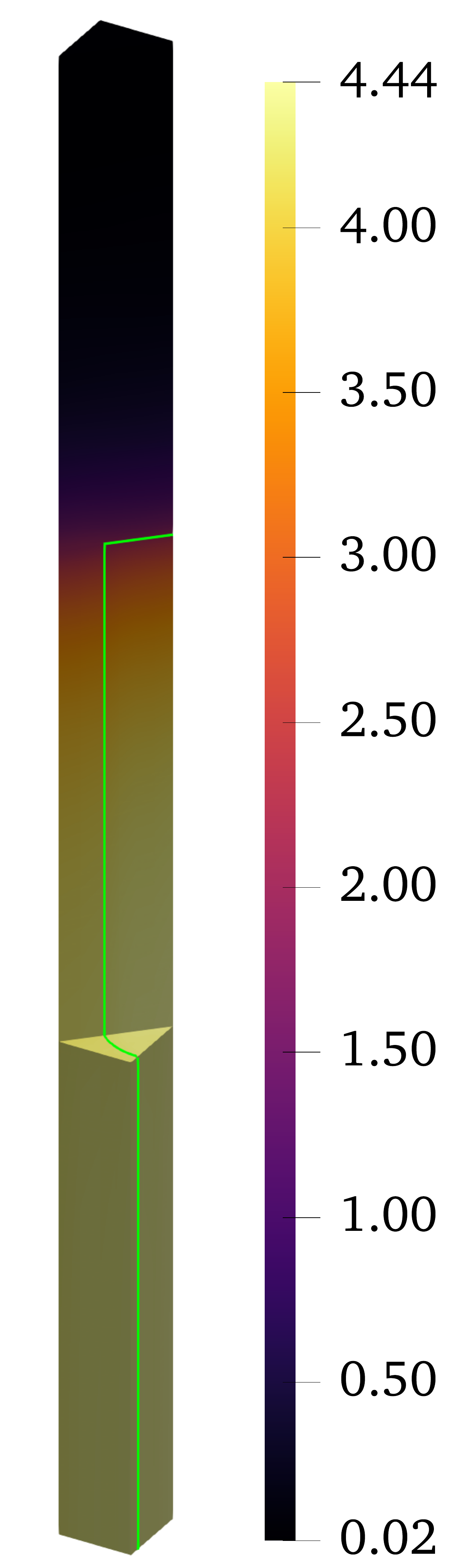}
		\caption{}
	\end{subfigure}
	\begin{subfigure}[t]{0.32\textwidth}
		\centering
		\includegraphics[width=0.425\textwidth]{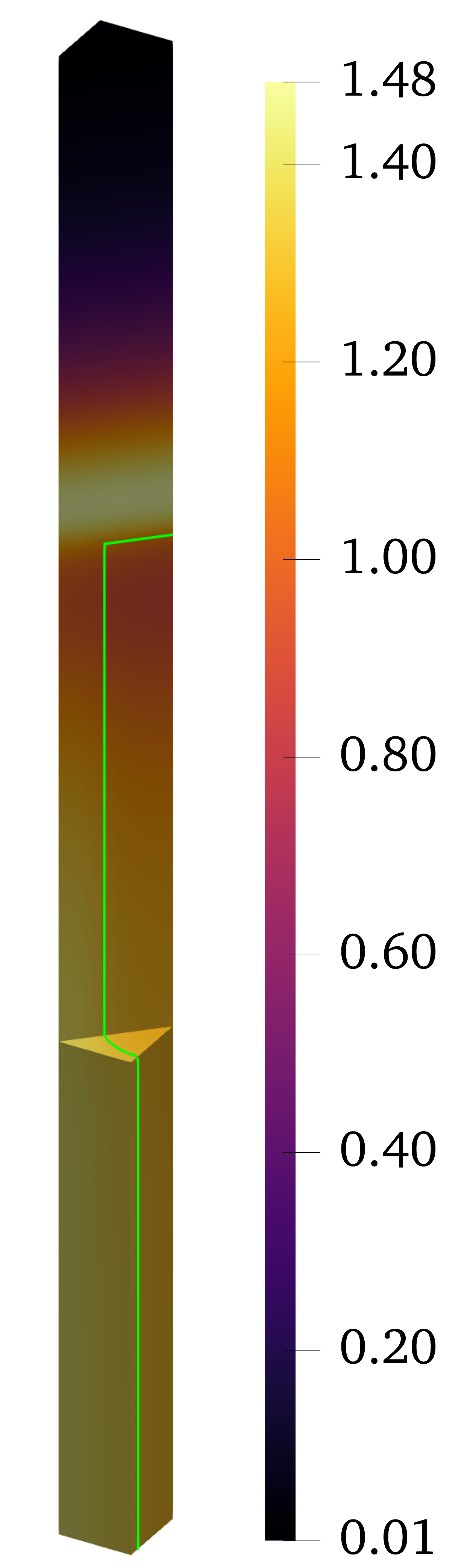}
		\caption{}
	\end{subfigure}
	\begin{subfigure}[t]{0.32\textwidth}
		\centering
		\includegraphics[width=0.425\textwidth]{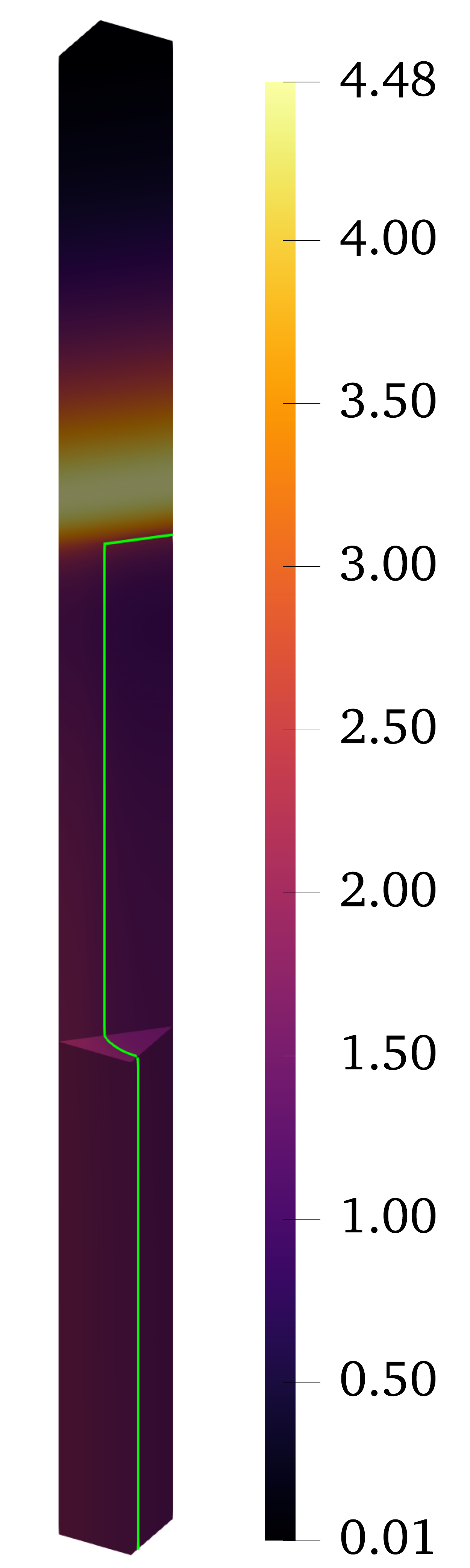}
		\caption{}
	\end{subfigure}
	\caption{Reference scalar fluxes for the nominal height 3D UO\textsubscript{2} C5G7 pin-cell in fast (a), epi-thermal (b), and thermal (c) energy groups, $g=1,\,6,\,7$. The cutaway reveals the radial distribution along the diagonal of the pin-cell, where the outline of the pin is drawn in green. For convenience, the $z$-dimension has been re-scaled (shrunk) by a ratio of 5:1.}
	\label{fig:pin-3D}
\end{figure}

Next, we consider a square pin-cell, fueled by UO\textsubscript{2} and moderated by light water, as defined in the (modified) C5G7 benchmark \cite{Lewis2005}. 
Moreover, to investigate how the pin height affects the quality of the low-rank decomposition, we grant that the nominal height of 42.84 cm may be adjusted by $\pm 50\%$ to achieve three test cases: namely, short (21.42 cm), medium (nominal), and tall (64.26 cm). In each case, the height of the axial reflector is set to the nominal value, 21.42 cm.
The radial mesh is shown in Figure \ref{fig:pin-2D} while the axial mesh is taken as a uniform grid of 0.1575 cm (equal to one-eighth of the 1.26 cm pitch) spacing. Meanwhile, the seven-group energy mesh is as prescribed by the benchmark, as are the corresponding cross sections. The 3D reference solution is shown in Figure \ref{fig:pin-3D} for the nominal height fuel pin.

\begin{figure}
	\centering
	\includegraphics[width=0.915\textwidth]{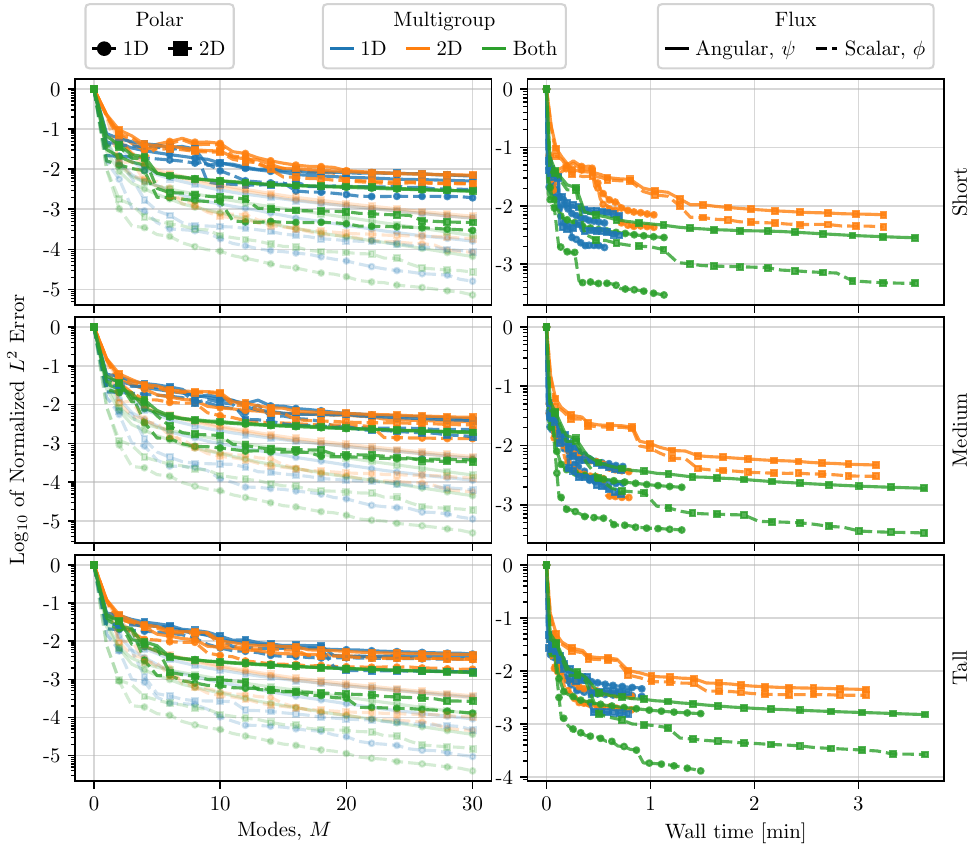}
	\caption{Convergence of 2D/1D Proper Generalized Decomposition (opaque lines) and Singular Value Decomposition (faint lines) for the C5G7 UO\textsubscript{2} pin-cell parameterized by pin height: short (21.42 cm), medium (42.84 cm), or tall (64.26 cm).}
	\label{fig:pin-2D1D}
\end{figure}

Plotting the SVD convergences in Figure \ref{fig:pin-2D1D}, we find again that---given 30 modes---the $L^2$ norm of the error in each decomposition is less than $10^{-3}$. Likewise, each PGD ROM converges within $10^{-2}$ of the 3D reference solution by the same metric. Moreover, while comparable, the fastest convergence with rank $M$ is observed in the case where both radial and axial submodels are multigroup; meanwhile, the difference between axial/-polar PGD is seen to be marginal. 
Regarding the convergence with runtime, we find---as in the Takeda LWR---that the axial PGD with 2D energy-dependence
is typically the least economical, while the axial-polar PGD where both submodels are multigroup is the most economical. 
The performance of other ROMs, meanwhile, is roughly comparable.
As to the pin height, its effect on each decomposition appears to be slight enough that the convergence with respect to either metric (rank or runtime) is comparable in each test case (short, medium, and tall pins). That said, there is a marginal improvement in the quality of the low-rank approximations (especially those of SVD) with increasing pin-height. Foreseeably, this is because an axially-homogeneous and -infinite problem with an axially-uniform source necessarily admits a rank-one solution (being axially constant)---in other words, lengthening the pin means there is proportionately less flux near the axial interface between fuel and reflector, where the distribution predictably requires more 2D/1D modes to resolve.

\section{Conclusion}
In summary, we have devised and implemented six PGD ROMs for 2D/1D decomposition. These are methodologically distinct, but practically analogous to, the conventional 2D/1D methods already popular in reactor physics. However, while the latter rest on physical
assumptions 
that 
limit 
the accuracy and reliability of the solution, the former assumes only that a sufficiently precise approximation can be achieved in a tractable number of modes $M$. Accordingly, 2D/1D PGD is hypothetically capable of converging 
to the 3D solution within an arbitrarily small precision (though in practice this convergence may be prohibitively slow), which is not the case of conventional 2D/1D methods. 

Moreover, beyond this key distinction, this
work introduces further novelties not found in previous 2D/1D methods. First, the latter assume that both the 2D and 1D components are functions of polar angle, which---if the 1D model is selected to be diffusion---is analogous to the axial PGD shown here.\footnote{More precisely, none of the PGD ROMs presented here grant that both the radial and axial modes depend on $\mu$. Nevertheless, the 1D---though not 1D($\mu$)---submodel of 2D($\mu$)/1D PGD is equivalent to a diffusion, or $P_1$, model, albeit one that can incorporate scattering of $L>1$, as explained in \ref{app:1D-P1}.}
However,
the axial-polar PGD demonstrates that $\mu$ can be separated from $\vect{\Omega}$ analogously to how 2D/1D methods separate $z$ from the 3D position vector.
Moreover, numerical experiments on the first Takeda benchmark (a small LWR) and a lattice of UO\textsubscript{2} pin-cells show this can
be practically advantageous, yielding markedly superior precision for an equal amount of wall time. Second, 2D/1D methods traditionally grant that both 2D and 1D equations are energy-dependent (multigroup). Presently, we consider the same, plus the two additional possibilities that $E$ is assigned uniquely to either the axial or radial domain. In this case, numerical performances on the two benchmarks
are comparable, with the ``group-wise'' decomposition often being the most economical (albeit by a small margin). Nevertheless, it is foreseeable that separating the energy dependence may be preferable for simulations with finer group-structures (the present examples having only two and seven groups respectively).

Given these unique advantages, one might expect 2D/1D PGD ROMs to compare favorably to existing methods. 
Alternatively, as these latter two features seem equally applicable to established 2D/1D methods---in that they apparently entail only the extension of ``transverse integration'' to include dimensions other than axial or radial space---future research could seek to incorporate these concepts into more traditional 2D/1D approaches. 
As either outcome stands to substantially advance the state-of-the-art in the  efficient modeling of 3D reactor physics,
we anticipate further investigation of 2D/1D PGD to be a fruitful direction of research, and view the present results as an encouraging first step in the development and application thereof.

\subsection{Possible Extensions}
Regarding 2D/1D PGD specifically, one clear opportunity for future work would be to enhance the performance
by the same means introduced by Nouy \cite{Nouy2010} and already demonstrated for neutron transport separated in energy by PGD in \cite{Dominesey2022,Dominesey2023}:
namely, Minimax PGD and the update step (foreseeably, here of the axial modes). Other options in this vein could include applying Subspace PGD \cite{Nouy2010} or the Arnoldi-like PGD of \cite{Tamellini2014}. Likewise, extending this ROM to $k$-eigenvalue problems---perhaps as in \cite{Dominesey2023}---would allow for wider applicability and more meaningful validation.
Beyond these algorithmic improvements,
the prospects  
for future research in reactor physics could include
incorporating prior knowledge of the solution  and implementing the performance-critical optimizations (especially preconditioning and parallelism) necessary to apply these ROMs to High Performance Computing applications and perform rigorous performance benchmarks.

Yet another direction could be to investigate whether features developed for conventional 2D/1D methods can be applied to enhance PGD. These could include 3D CMFD preconditioners, localization of the decomposition (perhaps to pin-cells, assemblies, and/or axial layers), and energy-dependent relaxation factors, as in  \cite{Collins2016}. Likewise, direct performance comparisons between the two approaches would be illuminating. That said, computing the requisite inner products for PGD may prove more complicated in the Method of Charcteristics, rather than finite element, discretization of space; as such, applying and modifying 2D/1D PGD to the former may constitute a separate research objective in and of itself. In any case, it is clear that---like the original 2D/1D ``fusion'' method of Cho et al. \cite{Cho2002}---the 2D/1D ROMs presented here stand to benefit considerably from the continued research of knowledgeable reactor physicists.

\section*{Acknowledgments}
This research was performed under appointment of the first author to the Rickover Fellowship Program in Nuclear Engineering sponsored by the Naval Reactors Division of the National Nuclear Security Administration.

\bibliographystyle{model1-num-names}
\bibliography{refs}

\appendix
\def\pull{\!\!\!}
\section{Relating the Axial and $P_1$ Equations}
\label{app:1D-P1}
The 1D submodel, Equation \ref{eq:prog-1D}, of the axial or 2D($\mu$)/1D PGD amounts to neutron transport with a two-point angular quadrature at $\mu=\pm 1$.
While this is straightforward to solve numerically, it
may be useful in some cases to recast this equation to an equivalent form resembling the familiar $P_1$ equations; this can be accomplished as follows.
For convenience, 
a more typical nomenclature is recovered 
by renaming $Z^\pm_{m,g}$ to the partial currents $J^\pm_{m,g}$, the even and odd axial fluxes $Z^{0}_{m,g}$ and $Z^{1}_{m,g}$ to the scalar flux and current $\phi_{m,g}$ and $J_{m,g}$, and the source term $\widetilde{q}^{(m^*)}_{\text{1D},g}$ to $Q_g^\pm$. Doing so, Equation \ref{eq:prog-1D} can be rewritten as
\begin{equation}
\label{eq:axial-submodel}
\begin{split}
\sum_{m=1}^M\Bigg[&\pm s_{\oneD,g}^{(m^*,m)}\frac{\partial}{\partial z} J_{m,g}^\pm(z)
+\widetilde{\Sigma}_{t,g}^{(m^*,m)}(z)J^\pm_{m,g}(z)
-\frac{1}{2}
\sum_{p=0}^1
(\pm 1)^{p}\sum_{g'=1}^G\widetilde{\Sigma}_{s,p,g'\rightarrow g}^{(m^*,m)}(z)\varphi_{m,g'}^p(z)\Bigg]
=Q^\pm_g(z)
\end{split}
\end{equation}
where
\begin{equation}
\varphi_{m,g}^p\equiv \begin{cases}
\phi_{m,g}(z)\equiv J^+_{m,g}(z)+J^-_{m,g}(z)\,, & \ell+k\text{ is even}\,,  \\
J_{m,g}(z)\equiv J_{m,g}^+(z)-J_{m,g}^-(z)\,,  & \ell+k\text{ is odd}\,.
\end{cases}
\end{equation}
Evidently, these two equations relating $J_{m,g}^+$ and $J_{m,g}^-$ (or, for brevity, Equation \ref{eq:axial-submodel}\textsuperscript{$\pm$})
can be transformed to those relating $\phi_{m,g}$ and $J_{m,g}$ by
adding Equations \ref{eq:axial-submodel}\textsuperscript{$+$} and \ref{eq:axial-submodel}\textsuperscript{$-$} to yield
\begin{equation}
\label{eq:axial-q0}
\sum_{m=1}^M\left[s_{\oneD,g}^{(m^*,m)}\frac{\partial }{\partial z}  J_{m,g}(z)+
\widetilde{\Sigma}_{t,g}
^{(m^*,m)}(z)\phi_{m,g}(z)-
\sum_{g'=1}^G
\widetilde{\Sigma}_{s,0,g'\rightarrow g}^{(m^*,m)}(z)\phi_{m,g'}(z)\right]
=Q_{g}^0(z)
\end{equation}
while likewise Equation \ref{eq:axial-submodel}\textsuperscript{$+$} minus Equation \ref{eq:axial-submodel}\textsuperscript{$-$} yields
\begin{equation}
\label{eq:axial-q1}
\sum_{m=1}^M\left[
s_{\oneD,g}^{(m^*,m)}\frac{\partial }{\partial z} \phi_{m,g}(z)+
\widetilde{\Sigma}_{t,g}
^{(m^*,m)}(z)J_{m,g}(z)-
\sum_{g'=1}^G
\widetilde{\Sigma}_{s,1,g'\rightarrow g}^{(m^*,m)}(z)J_{m,g'}(z)\right]
=Q_{g}^1(z)
\end{equation}
where $Q_g^0$ and $Q_g^1$ are the zeroth and first-order---synonymously, even and odd parity---moments of the source
\begin{equation}
\begin{split}
Q_g^0(z)&\equiv Q_g^+(z)+Q_g^-(z)\,, \\
Q_g^1(z)&\equiv Q_g^+(z)-Q_g^-(z)\,.
\end{split}
\end{equation}
Equations \ref{eq:axial-q0} and \ref{eq:axial-q1} can readily be recognized as a system of equations similar to the common $P_1$ equations in reactor physics. One could solve them in this form, or in the case of isotropic scattering ($L=0$) symbolically solve Equation \ref{eq:axial-q1} for $J_{m^*,g}$ and substitute the result into Equation \ref{eq:axial-q0}.
This reformulation of Equation \ref{eq:prog-1D} may be useful in future PGD implementations or 2D/1D methods, especially as the diffusion approximation is generally the preferred 1D model in the latter case.

\end{document}